\documentclass[10pt,french]{amsart}
\usepackage[francais]{babel}
\usepackage[utf8]{inputenc}
\usepackage{latexsym}
\usepackage[all]{xy}
\usepackage{amsmath}
\usepackage{amssymb}
\usepackage{amsfonts}

\usepackage{color}
\definecolor{niceblue}{rgb}{0,0,0.6}
\usepackage{latexsym}

\usepackage{url}
\usepackage[all]{xy}
\usepackage{longtable}
\usepackage{dsfont}
\usepackage{tikz}
\usetikzlibrary{decorations.pathmorphing}
\usetikzlibrary{calc}
\usetikzlibrary{decorations.markings}

\def\N{{\mathbb N}}
\def\Z{{\mathbb Z}}

\def\tf{$\sqcup\hskip-2.3mm\sqcap$}   
 
{\newtheorem{theorem}{Th\'eor\`eme}[section]   
   \newtheorem{proposition}[theorem]{Proposition}
   \newtheorem{lemma}[theorem]{Lemme}  
   \newtheorem{remark}[theorem]{Remarque} 
   \newtheorem{corollary}[theorem]{Corollaire} 

   \newtheorem{definition}[theorem]{D\'efinition}    
{\theoremstyle{definition}   
   }
{\theoremstyle{remark}

   \newtheorem*{preuve}{Preuve}

\author{Michel Vaqui\'e}
\address{Institut de Math\'ematiques de Toulouse UMR 5219, CNRS, Universit\'e de Toulouse,  
UPS, 118 route de Narbonne, F-31062 Toulouse Cedex 9, 
France} 
\email{vaquie@math.univ-toulouse.fr} 
\title[Module des diff\'erentielles]{Module des diff\'erentielles} 

   \thanks{$*$ Partially supported by the grant of the Agence Nationale de la Recherche “CatAG”ANR-17-CE40-0014.} 

\begin{document}

\begin{abstract} 
Pour tout objet $x$ dans une cat\'egorie $\mathbf C$ il est possible de d\'efinir la cat\'egorie des \emph{modules de Beck} au dessus de $x$, comme la cat\'egorie $Ab({\mathbf C}_{/x})$ des objets en groupe ab\'elien de la cat\'egorie ${\mathbf C }_{/x}$. 
Nous pouvons en d\'eduire, au moins pour toute cat\'egorie localement pr\'esentable, la notion de module \emph{cotangent} ou de module des \emph{diff\'erentielles} $\Omega _x$ de $x$ dans $Ab({\mathbf C}_{/x})$. 

Dans le cas de la cat\'egorie $\mathbf{Alg} _k$ des $k$-alg\`ebres commutatives au-dessus d'un anneau $k$, la cat\'egorie des modules de Beck $Ab({\mathbf{Alg}_k}_{/A})$ au-dessus d'une $k$-alg\`ebre $A$ est \'equivalente \`a la cat\'egorie $\mathbf{Mod}_A$ des $A$-modules et le module des diff\'erentielles est \'egal au module des \emph{diff\'erentielles de K\"ahler} de $A$. 

Le but de cet article est de montrer pour toute cat\'egorie localement pr\'esentable des r\'esultats qui g\'en\'eralisent les propri\'et\'es classiques des modules des diff\'erentielles de K\"ahler. 

\vskip .2cm

\noindent {\scshape Abstract}.
For any object $x$ in a category $\mathbf C$ it is possible to define the category of \emph{Beck modules} over $x$ as the category $Ab({\mathbf C}_{/x})$ of abelian group objects in the category ${\mathbf C}_{/x}$. 
We can deduce from this construction, at least for any locally presentable category, the notion of \emph{cotangent} module or module of \emph{differentials} $\Omega _x$ of $x$ in $Ab({\mathbf C}_{/x})$.  

In the case of the category $\mathbf{Alg} _k$ of commutative $k$-algebras over a ring $k$, the category of Beck modules $Ab({\mathbf{Alg}_k}_{/A})$ over a $k$-algebra $A$ is equivalent to the category $\mathbf{Mod}_A$ of $A$-modules and the cotangent module is equal to the module of \emph{K\"ahler differentials} of $A$. 

The aim of this article is to prove for any locally presentable category some results which generalize the classical properties of modules of K\"ahler differentials.  
 \end{abstract}


\vskip .2cm 
\date{Avril 2023}

\maketitle   

\tableofcontents

		      \section*{Introduction}
%

Soit $X$ une vari\'et\'e diff\'erentiable, analytique ou alg\'ebrique, pour tout point $x$ de $X$ nous pouvons d\'efinir un espace tangent $T _x X$ et pour tout morphisme $\xymatrix @C=8mm{g : X \ar[r] & Y}$ nous pouvons d\'efinir une application diff\'erentielle $\xymatrix @C=8mm{dg _x : T_xX \ar[r] & T_{g(x)}Y}$ qui est une application lin\'eaire qui approxime $g$ au voisinage de $x \in X$. 
Nous pouvons ainsi consid\'erer les notions d'espace tangent et de diff\'erentielle comme une  mani\`ere de lin\'eariser une cat\'egorie. 
Pour pouvoir g\'en\'eraliser cette construction nous allons d'abord rappeler la construction de l'espace tangent et de l'application diff\'erentielle dans le cadre alg\'ebrique.

Soient $k$ un anneau commutatif et $A$ une $k$-alg\`ebre, alors il existe un $A$-module $\Omega _{A/k}$, le \emph{module des diff\'erentielles de K\" ahler} ou \emph{module cotangent}, qui repr\'esente le foncteur 
$$\begin{array}{cccc}
Der _k (A,-) : & \mathbf{Mod}_A &  \longrightarrow & \mathbf{Ens} \\
 & M &\mapsto & Der_k(A,M),
\end{array}$$
o\`u pour tout $A$-module $M$ nous notons $Der_k(A,M)$ l'ensemble des \emph{$k$-d\'erivations} de $A$ \`a valeurs dans $M$.

Alors si $X$ est une vari\'et\'e alg\'ebrique sur $k$, nous pouvons d\'efinir le faisceau $\Omega _{X/k}$ des diff\'erentielles de K\" ahler de $X$ comme le faisceau quasi-coh\'erent sur $X$ d\'efini par le $A$-module $\Omega _{A/k}$ sur tout ouvert affine $Spec(A)$ de $X$, et le  fibr\'e tangent $\xymatrix @C=8mm{TX \ar[r] & X}$ est \'egal par d\'efinition au fibr\'e  $\mathbb{V}(\Omega _{X/k}) := \rm{Spec} \bigl ( Sym ( \Omega _{X/k} ) \bigr )$.

Pour tout morphisme $\xymatrix @C=8mm{f : A \ar[r] & B}$ dans la cat\'egorie $\mathbf{Alg} _k$ des $k$-alg\`ebres commutatives il existe une application de $B$-modules 
$$\xymatrix @C=8mm{ \delta _f : \Omega _A \otimes _A B \ar[r] & \Omega _B} \ ,$$ 
nous en d\'eduisons que pour tout morphisme $\xymatrix @C=8mm{g : X \ar[r] & Y}$ nous avons un diagramme commutatif 
$$\xymatrix @C=12mm @R=8mm{TX \ar[r] ^{dg} \ar[d]  & TY \ar[d] \\  
X \ar[r] ^g & Y
}$$ 
o\`u le morphisme $\xymatrix @C=8mm{TX \ar[r] & TY \times _Y X}$ correspond \`a l'application $\xymatrix @C=8mm{ g^*\bigl ( Sym(\Omega _Y) \bigr )\ar[r] & Sym(\Omega _X)}$ induite par l'application naturelle $\xymatrix @C=8mm{ \delta _g : g^*(\Omega _Y) \ar[r] & \Omega _X}$. 

\vskip .2cm

Le module cotangent, ou plus g\'en\'eralement le complexe cotangent qui est la version \emph{d\'eriv\'ee} du module cotangent dans le cas d'une vari\'et\'e singuli\`ere, permet d'\'etudier les d\'eformations de la vari\'et\'e, et plus pr\'ecis\'ement d'\'etudier les d\'eformations infinit\'esimales. 
Pour cela il est utile de consid\'erer les \emph{extensions de carr\'e nul} d'un anneau $A$. 
Nous rappelons que pour tout $A$-module $M$ nous pouvons munir le $A$-module $A \oplus M$ d'une structure de $k$-alg\`ebre commutative en posant $(a,m) . (a',m') = (aa',am'+a'm)$. 

Nous consid\'erons la $k$-alg\`ebre $B=A \oplus M$ comme un objet de la cat\'egorie ${\mathbf{Alg} _k} _{/A}$, c'est-\`a-dire comme une $k$-alg\`ebre munie d'un morphisme $\xymatrix @C=8mm{u_B : B \ar[r] & A}$, et dont les morphismes sont les triangles commutatifs 
$$\xymatrix @C=8mm @R=4mm{ 
B \ar[rr] ^f \ar[rd] _{u_B} && C \ar[ld]^{u_C} \\ 
& A
}$$ 

Nous remarquons alors que $B=A \oplus M$ n'est pas un objet arbitraire de la cat\'egorie ${\mathbf{Alg} _k}_{/A}$, mais est muni d'une structure \emph{d'objet en groupe ab\'elien}, et que le foncteur qui envoie $M$ sur $A \oplus M$ d\'efinit une \'equivalence de cat\'egories entre la cat\'egorie $\mathbf{Mod}_A$ des $A$-modules et la cat\'egorie $Ab({\mathbf{Alg} _k}_{/A})$ des objets en groupe ab\'elien dans ${\mathbf{Alg} _k}_{/A}$. 

Le foncteur oubli de la structure d'objet en groupe ab\'elien $\xymatrix @C=8mm{ U_A : Ab({\mathbf{Alg} _k}_{/A}) \ar[r] & {\mathbf{Alg} _k}_{/A}}$ admet un foncteur adjoint \`a gauche 
$\xymatrix @C=8mm{L_A : {\mathbf{Alg} _k}_{/A} \ar[r] & Ab({\mathbf{Alg} _k}_{/A})}$ 
et nous avons l'\'egalit\'e $\Omega _{A/k} = L_A( \star _{{\mathbf{Alg} _k}_{/A} })$, o\`u $\star _{{\mathbf{Alg} _k}_{/A} }$ est l'objet final 
$\xymatrix @C=8mm{id_A : A \ar[r] & A}$ de la cat\'egorie ${\mathbf{Alg} _k}_{/A}$. 

\vskip .2cm 

Il est possible de g\'en\'eraliser cette construction, nous pouvons d\'efinir les objets en groupe ab\'elien dans toute cat\'egorie $\mathbf C$ comme les objets $a$ de $\mathbf C$ tel que le pr\'efaisceau repr\'esentable associ\'e par le morphisme de Yoneda $\xymatrix @C=8mm{h_a : {\mathbf C} ^{op} \ar[r] & \mathbf{Ens}}$ est un pr\'efaisceau en groupes ab\'eliens. Alors pour tout objet $x$ de $\mathbf C$ nous d\'efinissons la cat\'egorie des \emph{modules de Beck} au-dessus de $x$ comme la cat\'egorie $Ab({\mathbf C}_{/x})$ des objets en groupe ab\'elien dans la cat\'egorie ${\mathbf C}_{/x}$. 
De plus la cat\'egorie $Ab({\mathbf C}_{/x})$ est une cat\'egorie additive, qui est ab\'elienne si la cat\'egorie $\mathbf C$ est une cat\'egorie exacte. 

Alors si le foncteur \emph{oubli} $\xymatrix @C=8mm{U_x : Ab({\mathbf C}_{/x}) \ar[r] & {\mathbf C}_{/x}}$ admet un adjoint \`a gauche $\xymatrix @C=8mm{L_x : {\mathbf C}_{/x} \ar[r] & Ab({\mathbf C}_{/x})}$, nous appelons \emph{module cotangent abstrait} ou \emph{module des diff\'erentielles} de $x$ l'image par $L_x$ de l'objet final $\star _{ {\mathbf C}_{/x}} = \xymatrix @C=8mm{ id_x : x \ar[r] & x}$ de ${\mathbf C}_{/x}$, et nous notons $\Omega _x = L_x (\star _{ {\mathbf C}_{/x}} )$. 

Nous renvoyons aux articles de J.M. Beck \cite{Be} et de D.G. Quillen \cite{Qu} pour le lien entre la cat\'egorie des objets en groupe ab\'elien et la cat\'egorie des modules et \`a l'article de M. Barr \cite{Ba1} pour les propri\'et\'es de la cat\'egorie des modules de Beck.   

\vskip .2cm 

Comme dans le cas des $k$-alg\`ebres commutatives, nous voulons alors d\'efinir pour tout morphisme $\xymatrix @C=8mm{f : x \ar[r] & y}$ dans $\mathbf C$ une application entre les modules des diff\'erentielles $\Omega _x$ et $\Omega _y$. 
Pour $f$ dans $\mathbf C$ il existe un foncteur additif   
$\xymatrix{f^* _{[Ab]}: Ab({\mathbf C}_{/y}) \ar[r] & Ab({\mathbf C}_{/x})}$ entre les cat\'egories des modules de Beck au dessus respectivement de $y$ et de $x$  et nous construisons une application naturelle 
$\xymatrix@C=8mm{ [{\delta} _f] : \Omega _x \ar[r] & f_{[Ab]}^* (\Omega _y)}$ dans $Ab({\mathbf C}_{/x})$. 
Alors si le foncteur $f^* _{[Ab]}$ admet un adjoint \`a gauche $\xymatrix{{f_!}_{[Ab]}: Ab({\mathbf C}_{/x}) \ar[r] & Ab({\mathbf C}_{/y})}$ nous en d\'eduisons une application 
$\xymatrix @C=8mm{ [\tilde{\delta _f}] :{f_!}_{[Ab]} ( \Omega _x) \ar[r] & \Omega _y}$ 
dans $Ab({\mathbf C}_{/y})$ qui correspond \`a l'application $\xymatrix @C=8mm{ \delta _f : \Omega _A \otimes _A B \ar[r] & \Omega _B}$ dans le cas de la cat\'egorie ${\mathbf{Alg} _k}$ des $k$-alg\`ebres commutatives. 

\vskip .2cm 

Dans le cadre de la cat\'egorie des $k$-alg\`ebres l'application $\delta _f$ donne des informations sur le morphisme $\xymatrix @C=8mm{f : A \ar[r] & B}$. 
En particulier nous avons deux suites exactes fondamentales entre les modules des diff\'erentielles, et un des objectif de cet article est de g\'en\'eraliser ces suites exactes au cadre d'une cat\'egorie $\mathbf C$ localement pr\'esentable quelconque. 

Pour \'enoncer ces r\'esultats nous devons d'abord d\'efinir  pour tout morphisme $\xymatrix @C=8mm{f : x \ar[r] & y}$ dans la cat\'egorie $\mathbf C$ l'analogue du $B$-module des diff\'erentielles relatives $\Omega _{B/A}$ d\'efini pour tout morphisme $\xymatrix @C=8mm{f : A \ar[r] & B}$ de $k$-alg\`ebres. 
Pour cela nous construisons le module $\Omega _f$ comme objet dans la cat\'egorie des modules de Beck au dessus de $f$ consid\'er\'e comme objet de la cat\'egorie ${\mathbf C}_{x/}$ (cf. d\'efinition \ref{def:module-des-differentielles-relatives}).  

\vskip .2cm 

Alors le th\'eor\`eme \ref{th:premiere-suite-exacte} g\'en\'eralise la \emph{premi\`ere suite exacte fondamentale} d\'efinie entre les modules des diff\'erentielles dans le cadre de ${\mathbf{Alg} _k}$,
le th\'eor\`eme \ref{th:deuxieme-suite-exacte} est inspir\'e par la \emph{deuxi\`eme suite exacte fondamentale} et la compatibilit\'e du module des diff\'erentielles avec la localisation, 
enfin le th\'eor\`eme \ref{th:troisieme-suite-exacte} est la g\'en\'eralisation du fait que les modules de diff\'erentielles commutent avec \emph{l'extension des scalaires}.

\vskip .2cm 

Dans la derni\`ere partie nous \'etudions comment ces r\'esultats s'\'enoncent dans les cas des cat\'egories des ensembles, des mono\"\i des commutatifs et nous revenons sur l'exemple de la cat\'egorie des $k$-alg\`ebres commutatives. 

\vskip .5cm

    
      \section{Rappels sur les foncteurs adjoints}    


Soient $\mathbf C$ et $\mathbf D$ deux cat\'egories, si nous avons une paire de foncteurs adjoints 
$$\xymatrix @C=12mm{
F : {\mathbf C} \ar@<1ex>[r] & {\mathbf D}:G \ar@<1ex>[l]_{\bot}}$$
nous notons respectivement $\xymatrix @C=8mm{\eta _x: x \ar[r] & GFx}$ et $\xymatrix @C=8mm{\epsilon _y : FGy \ar[r] & y}$ l'unit\'e et la counit\'e de l'adjonction. 
Nous rappelons que les morphismes compos\'es $\xymatrix @C=8mm{\epsilon _{Fx} \circ F \eta _x : Fx \ar[r] & Fx}$ et $\xymatrix @C=8mm{G \epsilon _y \circ \eta _{Gy} : Gy \ar[r] & Gy}$ sont les morphismes identit\'es. 

Les bijections naturelles entre les ensembles $Hom _{\mathbf C} (x, Gy)$ et $Hom _{\mathbf D} (Fx, y)$ sont d\'efinies par: 

$$\xymatrix @R=.1mm @C=2mm 
{Hom _{\mathbf C} (x, Gy) \ar@<1ex>[rr] && Hom _{\mathbf D} (Fx, y)\ar[ll] \\ 
\alpha \ar@{|->}[rr] && \epsilon_y \circ F\alpha \\ 
G\beta \circ \eta _x && \beta \ar@{|->}[ll]
}$$

\vskip .2cm 

Soit $\mathbf C$ une cat\'egorie, pour tout objet $x$ de $\mathbf C$ nous notons respectivement $\mathbf{C}_{/x}$ et  $\mathbf{C}_{x/}$ les cat\'egories des objets au-dessus de $x$, c'est-\`a-dire la cat\'egorie form\'ee des fl\`eches $\xymatrix @C=8mm{u:z \ar[r] & x}$, et la cat\'egorie des objets au-dessous de $x$, c'est-\`a-dire la cat\'egorie form\'ee des fl\`eches $\xymatrix @C=8mm{v:x \ar[r] & z}$.  

Soit $\xymatrix @C=8mm{f:x \ar[r] & y}$ un morphisme dans $\mathbf C$, alors si la cat\'egorie $\mathbf C$ admet des limites finies nous avons une paire de foncteurs adjoints 
$$\xymatrix{
f_! : {\mathbf C}_{/x} \ar@<1ex>[r] & {\mathbf C}_{/y} : f^* \ar[l]
}$$
d\'efinie par 
$$\xymatrix  @R=1mm @C=5mm {
f_! : (z \ar[r] ^-u & x) \ar@{|->}[rr] && ( z \ar[r]^-{f\circ u} & y) \\ 
f^* :  (t \ar[r] ^-v & y) \ar@{|->}[rr] && ( x\times _y t \ar[r]^-{p_1} & x)
}$$

De m\^eme si la cat\'egorie $\mathbf C$ admet des colimites finies nous avons une paire de foncteurs adjoints 
$$\xymatrix @C=12mm{
f_* : {\mathbf C}_{x/} \ar@<1ex>[r] & {\mathbf C}_{y/} : f^! \ar[l]
}$$
d\'efinie par 
$$\xymatrix  @R=1mm @C=5mm {
f_* : (x \ar[r] ^-v & t) \ar@{|->}[rr] && ( y \ar[r]^-{q_1} & y \coprod _x t) \\ 
f^! : (y \ar[r] ^-u & z) \ar@{|->}[rr] && ( x \ar[r]^-{u \circ f} & z)
}$$

En particulier dans les cas o\`u $y$ est l'objet final ou dans le cas o\`u $x$ est l'objet initial nous trouvons les paires de foncteurs adjoints suivantes: 

$$\xymatrix @C=12mm{
[x]_! : {\mathbf C}_{/x} \ar@<1ex>[r] & {\mathbf C} : [x]^*\ar[l]
}
\quad\quad\hbox{et}\quad\quad 
\xymatrix @C=12mm{
[y]_* : {\mathbf C} \ar@<1ex>[r] & {\mathbf C}_{y/} : [y]^! \ar[l]
}$$

\vskip .2cm 

Soit $\xymatrix @C=12mm{
F : {\mathbf C} \ar@<2pt>[r] & {\mathbf D}:G \ar@<1pt>[l]}$ 
une paire de foncteurs adjoints, alors pour tout $y$ dans $\mathbf D$ nous avons une nouvelle paire de foncteurs adjoints 
$\xymatrix @C=12mm{
F_{G(y)} : {\mathbf C}_{/G(y)} \ar@<2pt>[r] & {\mathbf D}_{/y}:G_y \ar@<1pt>[l]}$ 
d\'efinie par 
$$\xymatrix  @R=1mm @C=6mm {
F_{G(y)} : (x \ar[r] ^-v & G(y)) \ar@{|->}[rr] && ( F(x) \ar[r]^-{\epsilon _y \circ F(v)} & y ) \\ 
G_y : (z \ar[r] ^-u & y) \ar@{|->}[rr] && ( G(z) \ar[r]^-{G(u)} & G(y))\ .
}$$

De m\^eme pour tout $x$ dans $\mathbf C$ nous avons une paire de foncteurs adjoints 
$\xymatrix{
F_{x} : {\mathbf C}_{x/} \ar@<2pt>[r] & {\mathbf D}_{F(x)/}:G_{F(x)} \ar@<1pt>[l]}$ 
d\'efinie par 
$$\xymatrix  @R=1mm @C=6mm {
F_x : (x \ar[r] ^-v & z ) \ar@{|->}[rr] && ( F(x) \ar[r]^-{F(v)} & F(z) ) \\ 
G_{F(x)} : (F(x) \ar[r] ^-u & y) \ar@{|->}[rr] && ( x \ar[r]^-{G(u) \circ \eta _x} & G(y)) \ .
}$$

\vskip .5cm 

\begin{remark}\label{rmq:foncteurs-induits} 
Consid\'erons un morphisme $\xymatrix @C=8mm{g:z \ar[r] & y}$ dans une cat\'egorie $\mathbf C$ admettant des colimites finies et la paire de foncteurs adjoints 
$\xymatrix @C=12mm{ F = [z] _* : {\mathbf C} \ar@<1ex>[r] & {\mathbf E} = {\mathbf C}_{z/} : [z]^! = G  \ar[l]}$. 
L'image de $(\xymatrix @C=8mm{z \ar[r] ^g & y})$ par le foncteur $G$ est $y$ et nous en d\'eduisons la nouvelle paire de foncteurs adjoints 
$\xymatrix @C=12mm{ F_y : {\mathbf C} _{/y} \ar@<1ex>[r] & {\mathbf E}_{/g} : G_g  \ar[l]}$
d\'efinie de la mani\`ere suivante: 
$$\xymatrix @R=2mm @C=1mm{ 
&&&&&&&&&& z \ar[dd] ^{q_1}  \ar@{=}[rrr] &&& z \ar[dd] ^{g}  &&&&&& z \ar[dd] ^{u}  \ar@{=}[rrr] &&& z \ar[dd] ^{g}  \\ 
F_y : & ( b \ar[rrr] ^{l} &&& y ) & \ar@{|->}[rrrr]  &&&&&&&&&&& { }   \hbox{et} && G_g : &&&&&   &  \ar@{|->}[rrrr] &&&& & ( a \ar[rrr] ^{h} &&& y )  \\ 
&&&&&&&&&& z \coprod b \ar[rrr] _{(g,l)} &&& y  &&&&&& a\ar[rrr] _{a} &&& y}
$$

La cat\'egorie ${\mathbf E}_{/g}$ d\'efinie par ${\mathbf E}_{/g} = {\big( {\mathbf C} _{z/} \big)} _{/g}$ est \'egale \`a la cat\'egorie ${\big( {\mathbf C} _{/y} \big)} _{g/}$, c'est la cat\'egorie des diagrammes $\xymatrix @C=8mm{z \ar[r] ^{v_a} & a \ar[r] ^{u_a} & y}$ avec $u_a \circ v_a = g$. 

La paire de foncteurs adjoints 
$\xymatrix @C=12mm{ F_y : {\mathbf C} _{/y} \ar@<1ex>[r] & {\mathbf E}_{/g} : G_g  \ar[l]}$ 
d\'efinie pr\'ec\'edemment peut aussi \^etre d\'efinie comme la paire 
$\xymatrix @C=12mm{ F_y = [g] _* : {\mathbf C} _{/y} \ar@<1ex>[r] & {\mathbf E}_{/g} : = [g] ^! = G_g  \ar[l]}$ . 
\end{remark}

\vskip .5cm


\section{Modules de Beck} 


\subsection{Objets en groupe ab\'elien et modules de Beck} 

Pour toute cat\'egorie $\mathbf Y$ la cat\'egorie $Ab(\mathbf{Y})$ des objets en groupe ab\'elien dans $\mathbf Y$ est la sous-cat\'egorie form\'ee des objets $y$ de $\mathbf Y$ tels que le pr\'efaisceau $\xymatrix @C=8mm{h_y: \mathbf{Y}^{op} \ar[r] & \mathbf{Ens}}$ repr\'esent\'e par $y$ est un pr\'efaisceau en groupes ab\'eliens. 
Cela signifie que pour tout $z$ dans $\mathbf Y$ l'ensemble $h_y(z)= Hom_{\mathbf{Y}}(z,y)$ est un groupe ab\'elien et que pour tout morphisme $\xymatrix @C=8mm{u:z \ar[r] & z'}$ dans $\mathbf Y$ l'application associ\'ee $\xymatrix @C=8mm{u^*: h_y(z') \ar[r] & h_y(z)}$ est un morphisme de groupes. 
Un morphisme dans $Ab(\mathbf{Y})$ est un morphisme de pr\'efaisceaux en groupes. 

Si la cat\'egorie $\mathbf Y$ admet des limites finies, il suffit de supposer que $\mathbf Y$ admet des produits finis et a un objet final $\star _{\mathbf Y}$, un objet en groupe ab\'elien est un objet $y$ de $\mathbf Y$ et la donn\'ee de morphismes $\xymatrix @C=8mm{m_y: y\times y \ar[r] & y}$, $\xymatrix @C=8mm{i_y: y \ar[r] & y}$ et $\xymatrix @C=8mm{e_y : \star _{\mathbf Y} \ar[r] & y}$ tels que $m_y$ d\'efinit une loi de groupes ab\'eliens avec $i_y$ pour l'inverse et $e_y$ pour l'\'el\'ement neutre. 
Un morphisme dans $Ab(\mathbf{Y})$ de $y$ dans $y'$ est un morphisme $\xymatrix @C=8mm{f:y \ar[r] & y'}$ dans $\mathbf Y$ tel que nous ayons 
$f \circ m_y = m_{y'} \circ (f \times f)$, $f \circ i_y = i_{y'} \circ f$ et $f \circ e_y = e_{y'}$.

Nous rappelons le r\'esultat suivant (cf. \cite{Ba1}, Th\'eor\`emes 1.5 et 2.4 du chapitre 2).

\begin{theorem} 
Pour toute cat\'egorie $\mathbf Y$ ayant des limites finies la sous-cat\'egorie $Ab({\mathbf Y})$ des objets en groupe ab\'elien est une cat\'egorie additive, qui est localement pr\'esentable si $\mathbf Y$ l'est. 

De plus si la cat\'egorie $\mathbf Y$ est exacte, la cat\'egorie $Ab({\mathbf Y})$ est une cat\'egorie ab\'elienne. 
\end{theorem} 

Soit $\mathbf Y$ une cat\'egorie admettant des limites finies, alors l'\emph{\'el\'ement nul} $0_{\mathbf Y}$ de la cat\'egorie additive $Ab({\mathbf Y})$ est l'objet final $\star _{\mathbf Y}$ muni de la structure de groupe triviale. 

Pour tout $x$ dans $Ab({\mathbf Y})$, le morphisme $\xymatrix @C=8mm{x \ar[r] & 0_{\mathbf Y}}$ dans $Ab({\mathbf Y})$ est le morphisme induit par le morphisme $\xymatrix @C=8mm{x \ar[r] & \star _{\mathbf Y}}$ dans $\mathbf Y$, et pour tout $y$ dans $Ab({\mathbf Y})$ le morphisme $\xymatrix @C=8mm{0_{\mathbf Y} \ar[r] & y}$ dans $Ab({\mathbf Y})$ est le morphisme induit par l'\'el\'ement neutre $\xymatrix @C=8mm{e_y: \star _{\mathbf Y} \ar[r] & y}$ de la structure de groupe de $y$. 

En particulier l'\'el\'ement neutre $0$ du groupe ab\'elien $Hom _{Ab({\mathbf Y})}(x,y)$ est induit par un morphisme $g$ dans $\mathbf Y$ qui se factorise par l'objet final $\star _{\mathbf Y}$: 
$\xymatrix{ x \ar[r] \ar@/^/[rr]^g & \star _{\mathbf Y} \ar[r] &y}$.

\vskip .2cm 

Il existe un foncteur \emph{oubli} $\xymatrix @C=8mm{U: Ab({\mathbf Y}) \ar[r] & {\mathbf Y}}$, et si la cat\'egorie ${\mathbf Y}$ est localement pr\'esentable le foncteur $U$ admet un adjoint \`a gauche $\xymatrix @C=8mm{L:{\mathbf Y} \ar[r] & Ab({\mathbf Y})}$, le foncteur \emph{objet en groupe ab\'elien libre} (cf. \cite{Ba1}, Exercice 1.5-3 du chapitre 6). 

Un foncteur exact \`a gauche $\xymatrix @C=8mm{G:{\mathbf Y} \ar[r] & {\mathbf X}}$ pr\'eserve les limites finies et l'image d'un objet en groupe ab\'elien est un objet en groupe ab\'elien, il induit alors un foncteur additif $\xymatrix @C=8mm{G_{[Ab]}: Ab({\mathbf Y}) \ar[r] & Ab({\mathbf X})}$ 
tel que le diagramme suivant est commutatif 
$$\xymatrix @C=12mm{ 
Ab(\mathbf{Y}) \ar[r] ^-{{G}_{[Ab]}} \ar[d] ^{U} & Ab(\mathbf{X}) \ar[d] ^{U} \\ 
 \mathbf{Y} \ar[r] ^-{G} & \mathbf{X} 
}$$

\vskip .5cm 

Soit $\mathbf C$ une cat\'egorie admettant des limites finies. 

\begin{definition} 
La cat\'egorie des modules de Beck au-dessus d'un objet $x$ de $\mathbf C$ est la cat\'egorie des objets en groupe ab\'elien dans la cat\'egorie $\mathbf{C}_{/x}$ des objets au-dessus de $x$: $B_{\mathbf C}(x) = Ab({\mathbf C}_{/x})$. 
\end{definition} 

\vskip .2cm 

Nous pouvons d\'efinir un foncteur \emph{oubli} $U_x$ de la cat\'egorie des modules de Beck $Ab(\mathbf{C}_{/x})$ au-dessus de $x$ dans la cat\'egorie $\mathbf{C}_{/x}$. 
Si la cat\'egorie $\mathbf C$ est localement pr\'esentable, il en est de m\^eme de la cat\'egorie $\mathbf{C}_{/x}$ et le foncteur $U_x$ admet un adjoint \`a gauche $L_x$, qui est appel\'e \emph{foncteur d'ab\'elianisation} et est not\'e $Ab_x$ dans \cite{Fr}: 
$$\xymatrix @C=12mm{
L_x : {\mathbf C}_{/x} \ar@<1ex>[r] & Ab({\mathbf C}_{/x}):U_x \ar@<1ex>[l]-_{\bot}}$$

\vskip .5cm

Soit $\xymatrix @C=8mm{G:\mathbf{D} \ar[r] & \mathbf{C}}$ un foncteur exact \`a gauche, pour tout $y$ dans $\mathbf D$ le foncteur $\xymatrix @C=8mm{G_y:\mathbf{D}_{/y} \ar[r] & \mathbf{C}_{/G(y)}}$ pr\'eserve l'objet final et les produits finis, l'image d'un objet en groupe ab\'elien est un objet en groupe ab\'elien. 
Nous avons donc un foncteur ${G_y}_{[Ab]}$ de la cat\'egorie des modules de Beck au-dessus de $y$, $B_{\mathbf D}(y) =Ab(\mathbf{D}_{/y})$, dans  la cat\'egorie des modules de Beck au-dessus de $G(y)$, $B_{\mathbf C}(G(y)) =Ab(\mathbf{C}_{/G(y)})$, qui est additif,  tel que le diagramme suivant soit commutatif: 

$$\xymatrix @C=12mm{ 
Ab(\mathbf{D}_{/y}) \ar[r] ^-{{G_y}_{[Ab]}} \ar[d] ^{U_y} & Ab(\mathbf{C}_{/G(y)}) \ar[d] ^{U_{G(y)}} \\ 
 \mathbf{D}_{/y} \ar[r] ^-{G_y} & \mathbf{C}_{/G(y)}
}$$

\vskip .5cm 

Nous noterons $[a]$ un \'el\'ement de la cat\'egorie $Ab({\mathbf C}_{/x})$, cela correspond \`a un objet $a$ de $\mathbf C$ avec un morphisme $\xymatrix @C=8mm{u_a:a \ar[r] & x}$, et \`a la donn\'ee des morphismes $\xymatrix @C=8mm{m_a: a \times _x a \ar[r] & a}$, $\xymatrix @C=8mm{i_a:a \ar[r] & a}$ et $\xymatrix @C=8mm{e_a: x \ar[r] & a}$ donn\'es par la structure de groupe.  
Nous avons l'\'egalit\'e $u_a \circ e_a = id  _x$, la cat\'egorie $Ab({\mathbf C}_{/x})$ est une sous-cat\'egorie de la cat\'egorie ${\mathbf C}_{x/ . /x}$ form\'ee des diagrammes commutatifs $\xymatrix{ x \ar[r] \ar@/^/[rr]^{id _x} & a \ar[r] &x}$, 
nous notons $[a] = (\xymatrix @C=10mm{a \ar[r] _-{u_a} & x \ar@{.>}@/_/[l]_-{e_a}})$.  

L'objet nul $0_{{\mathbf C}_{/x}}$ de la cat\'egorie additive $Ab({\mathbf C}_{/x})$ correspond \`a 
$[x] = (\xymatrix @C=10mm{x \ar[r] _-{id_x} & x \ar@{.>}@/_/[l]_-{id_x}})$.  

De m\^eme nous noterons $[{\alpha}]$ un morphisme de $Hom _{Ab({\mathbf C}_{/x})}([a] , [b] )$, c'est la donn\'ee d'un morphisme $\alpha$ dans $Hom _{\mathbf C} (a,b)$ tel que nous ayons le diagramme commutatif 
$$
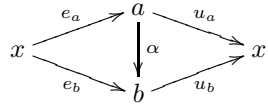
\xymatrix @C=12mm @R=2mm{ & a \ar[dd]^-{\alpha} \ar[rd]^-{u_a} \\ 
x \ar[ru]^-{e_a} \ar[rd]_-{e _b} && x \\ 
& b \ar[ru]_-{u_b}}$$
et compatible aux morphismes $m_a$ et $m_b$, $i_a$ et $i_b$.

\begin{lemma}\label{lem:morphisme-nul-0}  
Soit $\alpha$ un morphisme dans $Hom _{Ab({\mathbf C}_{/x})}([a] , [b] )$, alors avec les notations pr\'ec\'edentes nous avons 
$$[{\alpha}] = [0] \ \Longleftrightarrow \ \alpha = e_b \circ u_a \ .$$  
\end{lemma}  

\begin{preuve}. 
Le morphisme $[{\alpha}]$ est le morphisme nul dans le groupe ab\'elien $Hom _{Ab({\mathbf C}_{/x})}([a],[b])$ s'il se factorise par $0_{{\mathbf C}_{/x}}$, c'est-\`a-dire si nous pouvons trouver un diagramme commutatif  
$$\xymatrix @C=18mm @R=6mm{ & a \ar[d]^-{\alpha _1} \ar[rd]^-{u_a} \\ 
x \ar[ru]^-{e_a} \ar[rd]_-{e _b} \ar[r]^(.6){id _x} & x \ar[r]^(.4){id _x} \ar[d]^-{\alpha _2} & x \\ 
& b \ar[ru]_-{u_b}}$$

\hfill \tf 
\end{preuve} 

\vskip .5cm

\begin{definition} 
\begin{enumerate} 
\item Le \emph{module des diff\'erentielles}, ou \emph{module cotangent abstrait}, $\Omega _x$ de $x$ est le module de Beck au-dessus de $x$ obtenu comme image par $L_x$ de l'objet final $\star _{{\mathbf C}_{/x}} = (\xymatrix @C=8mm {x \ar[r] ^-{id_x} & x})$ de $\mathbf{C}_{/x}$, 
$\Omega _x = L_x ( \star _{{\mathbf C}_{/x}})$.  
\item Pour tout $[b]$ appartenant \`a $Ab({\mathbf C}_{/x})$, nous appelons ensemble des \emph{${\mathbf C}$-d\'erivations de $x$ \`a valeurs dans $[b]$} l'ensemble 
$$Der_{\mathbf C}(x,[b]):= Hom_{Ab({\mathbf C}_{/x})} (\Omega _x ,[b]) \simeq Hom_{{\mathbf C}_{/x}} (\star _{{\mathbf C}_{/x}}, b)$$  
avec $b=U_x([b])$. 
\end{enumerate} 
\end{definition} 

Le module des diff\'erentiel $\Omega _x$ correspond \`a un \'el\'ement 
$[{\omega _x}] = (\xymatrix @C=10mm{\omega _x \ar[r] _-{u_{\omega _x}} & x \ar@{.>}@/_/[l]_-{e_{\omega _x}}})$, avec $\omega _x$ dans $\mathbf C$, et par d\'efinition nous avons 
$U_x ( \Omega _x) = (\xymatrix @C=10mm{\omega _x \ar[r] ^-{u_{\omega _x}} & x})$. 
En particulier l'unit\'e $\eta _{\star _{{\mathbf C}_{/x}}}$ de l'adjonction $L_x \dashv U_x$ correspond \`a un morphisme $\xymatrix @C=8mm{ \eta _x :x \ar[r] & \omega _x}$ v\'erifiant $u_{\omega _x} \circ \eta _x = id _x$. 

\vskip .2cm 

Une d\'erivation $\theta$ dans $Der_{\mathbf C}(x,[b])$ correspond \`a un morphisme $[\delta]$ dans $Hom_{Ab({\mathbf C}_{/x})} ([\omega _x] ,[b])$, c'est-\`a-dire correspond \`a un morphisme $\xymatrix @C=8mm{ \delta : \omega _x \ar[r] & b}$ tel que le diagramme suivant soit commutatif 
$$\xymatrix @C=12mm @R=2mm{ & \omega _x \ar[dd]^-{\delta} \ar[rd]^-{u_{\omega _x}} \\ 
x \ar[ru]^-{e_{\omega _x}} \ar[rd]_-{e _b} && x \\ 
& b \ar[ru]_-{u_b}}$$

Par adjonction le morphisme $[\delta]$ correspond \`a un morphisme $\beta$  dans $Hom_{{\mathbf C}_{/x}} (\star _{{\mathbf C}_{/x}}, b)$, qui est d\'efini par $\beta = U([\delta ]) \circ \eta _{\star _{{\mathbf C}_{/x}}}$, c'est-\`a-dire \`a un morphisme $\xymatrix @C=8mm{ \beta :x \ar[r] & b}$ v\'erifiant $u_b \circ \beta = id _x$ tel que nous ayons l'\'egalit\'e $\beta = \delta \circ \eta _x$.  
Et r\'eciproquement par adjonction nous associons \`a un morphisme $\beta$ dans $Hom_{{\mathbf C}_{/x}} (\star _{{\mathbf C}_{/x}}, b)$ le morphisme $[\delta ]$ dans $Hom_{Ab({\mathbf C}_{/x})} ([\omega _x] ,[b])$ d\'efini par $[\delta ] = \epsilon_{[b]} \circ L_x(\beta )$.

\begin{proposition}\label{prop:morphisme-nul}  
Un \'el\'ement $\theta$ du groupe ab\'elien $Der_{\mathbf C}(x,[b])$ est l'\'el\'ement nul si et seulement si le morphisme $\beta$ de $Hom_{{\mathbf C}_{/x}} (\star _{{\mathbf C}_{/x}}, b)$ correspondant \`a $\theta$ est \'egal au morphisme $e_b$ \'el\'ement neutre de la structure de groupe ab\'elien de $[b]$. 
\end{proposition}  

\begin{preuve} 
Par d\'efinition l'\'el\'ement $\theta$ est l'\'el\'ement nul si nous avons $[\delta ] = [0]$.  

D'apr\`es le lemme \ref{lem:morphisme-nul-0}, si $[\delta ] = [0]$ nous avons $\delta = e_b \circ u_{\omega _x}$ et nous en d\'eduisons $\beta = e_b \circ u_{\omega _x} \circ \eta _x = e_b$. 

R\'eciproquement si $\beta = e_b$ alors $\beta$ est l'image par $U_x$ du morphisme nul $\xymatrix @C=8mm{[0] : 0_{{\mathbf C}_{/x}} \ar[r] & [b]}$, nous avons alors le diagramme commutatif suivant  
$$\xymatrix@C=12mm{ 
\omega _x = L_xU_x 0_{{\mathbf C}_{/x}} \ar[r]^-{L_x(\beta )} \ar[d] _-{\epsilon _{0_{{\mathbf C}_{/x}}}} \ar[rd] ^-{[\delta ]} & L_xU_x [b] \ar[d] ^-{\epsilon _{[b]}} \\ 
0_{{\mathbf C}_{/x}} \ar[r] _-{[0]} & [b] 
}$$
et nous en d\'eduisons que $[\delta ] = [0]$. 

\hfill \tf 
\end{preuve} 

\vskip .2cm 

\begin{remark} 
Si $x$ est l'objet initial de la cat\'egorie $\mathbf C$ pour tout objet $\xymatrix @C=8mm{u_b : b \ar[r] & x}$ dans ${\mathbf C}_{/x}$ il existe un seul morphisme $\xymatrix @C=8mm{\beta : x \ar[r] & b}$ tel que $u_b \circ \beta = id_x$, par cons\'equent le groupe ab\'elien $Der_{\mathbf C}(x,[b])$ est r\'eduit \`a l'\'el\'ement nul. 
En particulier nous avons $\Omega _x = (0)$. 
\end{remark} 

\vskip .5cm

\subsection{Cat\'egorie en groupe}

Nous consid\'erons maintenant une cat\'egorie $\mathbf{C}$ telle que tout objet $a$ admet une structure de groupe, non n\'ecessairement commutatif, et telle que les morphismes respectent cette structure. 
Plus pr\'ecis\'ement nous supposons que pour tout objet $a$ dans $\mathbf{C}$ il existe des morphismes $\xymatrix @C=8mm{\mu _a : a \times a \ar[r] & a}$, $\xymatrix @C=8mm{\iota _a : a \ar[r] & a}$ et $\xymatrix @C=8mm{\varepsilon _a : \star _{\mathbf C} \ar[r] & a}$ d\'efinissant une structure d'objet en groupe et que pour tout morphisme $\xymatrix @C=8mm{f : b \ar[r] & a}$ dans $\mathbf C$ nous avons $f \circ \mu _b = \mu _a \circ ( f \times f)$, $f \circ \iota _b = \iota _a \circ f$ et $f \circ \varepsilon _b = \varepsilon _a$.  
Et nous notons $\xymatrix @C=8mm{c_a : a \times a \ar[r] & a \times a}$ le morphisme qui \'echange les facteurs, en particulier la loi de groupe est commutative si nous avons $\mu _a \circ c_a = \mu _a$. 

\vskip .2cm 

Soit $x$ un objet dans $\mathbf C$ et nous notons comme pr\'ec\'edemment $m_a$, $i_a$ et $e_a$ les morphismes d\'efinissant la structure d'objet en groupe ab\'elien d'un objet $[a]$ dans $Ab({\mathbf C}_{/x})$, avec $\xymatrix @C=8mm{u _a : a \ar[r] & x}$.  

Nous avons le diagramme commutatif 
$$\xymatrix @C=12mm{ 
a \times _x a \ar[r] ^-{p_1} \ar[d] _{p_2} & a \ar[d] ^{u_a} \\ 
a \ar[r] ^{u_a} & x 
}$$  
et nous notons $\xymatrix @C=8mm{\mu _a ^{(3)} : a \times a \times a \ar[r] & a}$ le morphisme $\mu _a ( \mu _a , id _a )  =  \mu _a ( id _a , \mu _a )$, $\xymatrix @C=8mm{ g_a : a \ar[r] & a}$ le morphisme $g_a = e_a \circ u_a$,  et $\xymatrix @C=8mm{ q _a : a \ar[r] & a}$ le morphisme $g_a \circ p_1 = g_a \circ p_2$

\begin{proposition} \label{prop:loi-de-groupe}
Nous avons l'\'egalit\'e 
$m_a = \mu _a ^{(3)} \circ \phi$ 
o\`u le morphisme $\xymatrix @C=8mm{\phi : a \times _x a \ar[r] & a \times a \times a}$ est d\'efini par $\phi = (p_2, \iota _a \circ q _a , p_1)$. 
\end{proposition} 

\begin{preuve} 
C'est un r\'esultat classique qui est une cons\'equence des remarques suivantes. 

Comme les morphismes dans $\mathbf C$ respectent la loi de composition $\mu _a$ nous avons la commutativit\'e du diagramme suivant 
$$\xymatrix @C=5mm @R=3mm{
(a\times _x a ) \times ( a\times _x a) \ar[dd] _{m_a \times m_a} \ar[rr] ^-{\simeq} && (a \times a ) \times _{(x \times x)} (a \times a) \ar[dd] ^{\mu _a \times \mu _a} \\ \\ 
a \times a \ar[rd] ^{\mu _a} && a \times _x a \ar[ld] _{m_a} \\ 
& a
}$$
et nous en d\'eduisons  que le morphsime compos\'e $\xymatrix @C=7mm{
a\times  a \ar[r] ^-{\psi} & (a \times a ) \times _{(x \times x)} (a \times a) \ar[rrr] ^-{m_a \circ (\mu _a \times \mu _a)}
&&& a}$, 
avec $\psi = ((g_a \circ p_2, p_1), (p_2,g_a \circ p_1))$, est \'egal \`a $\mu _a \circ c_a$. 

\hfill \tf 
\end{preuve} 

\vskip .2cm 

Si nous supposons que les objets de la cat\'egorie $\mathbf C$ sont des ensembles munis d'une loi de groupe not\'ee multiplicativement les relations pr\'ec\'edentes peuvent s'\'ecrire: 
$m_a(r_1 . s_1 , r_2 . s_2) = m_a(r_1,s_1) . m_a(r_2,s_2)$ 
pour $r_1,r_2,s_1,s_2$ dans l'ensemble $a$ avec $u_a(r_1)=u_a(s_1)$ et $u_a(r_2)=u_a(s_2)$, et 
$m_a(g_a(s).r , s.g_a(r)) = s.r$ pour $r$ et $s$ quelconques dans $a$, et l'\'egalit\'e de la proposition peut s'\'ecrire 

$$m_a(r,s) = s . e_a(t)^{-1} . r \ ,$$
pour $r$ et $s$ appartenant \`a l'ensemble $a$ avec $t= u_a(r) =u_a(s)$. 

\vskip .2cm 

\begin{corollary} \label{cor:sous-categorie-pleine} 
La cat\'egorie $Ab({\mathbf C}_{/x})$ est une sous-cat\'egorie pleine de la cat\'egorie ${\mathbf C}_{x/ . /x}$ form\'ee des diagrammes commutatifs $\xymatrix @C=8mm{ x \ar[r] ^{e_a} & a \ar[r] ^{u_a} &x}$, avec $u_a \circ e_a = id _x$.   
\end{corollary} 

\begin{preuve} 
En effet nous d\'eduisons de la proposition pr\'ec\'edente que la structure d'objet en groupe ab\'elien de $[a]$ est enti\`erement d\'etermin\'ee par les morphisme $\xymatrix @C=8mm{e_a : x \ar[r] & a}$ et $\xymatrix @C=8mm{u_a : a \ar[r] & x}$. 

\hfill \tf 
\end{preuve} 

\vskip .5cm

\subsection{Suite exacte des modules des diff\'erentielles} 

Soit $\mathbf C$ une cat\'egorie localement pr\'esentable, et soit $\xymatrix@C=8mm{f:x \ar[r] & y}$ dans $\mathbf C$, nous voulons comparer les modules des diff\'erentielles $\Omega _x$ et $\Omega _y$ et d\'efinir un module des diff\'erentielles relatives $\Omega _{x/y} = \Omega _f$.  

Le foncteur $\xymatrix{f^*: \mathbf{C}_{/y} \ar[r] & \mathbf{C}_{/x}}$ est un adjoint \`a droite, il est donc exact \`a gauche et induit un foncteur additif   
$\xymatrix{f^* _{[Ab]}: Ab({\mathbf C}_{/y}) \ar[r] & Ab({\mathbf C}_{/x})}$ 
tel que le diagramme suivant soit commutatif 
$$\xymatrix @C=12mm{ 
Ab(\mathbf{C}_{/y}) \ar[r] ^-{f_{[Ab]}^*} \ar[d] ^{U_y} & Ab(\mathbf{C}_{/x}) \ar[d] ^{U_x} \\ 
 \mathbf{C}_{/y} \ar[r] ^-{f^*} & \mathbf{C}_{/x}
}$$

\vskip .2cm 

Nous avons un morphisme canonique $\xymatrix @C=8mm{\vartheta _f : f_!(\star _{{\mathbf C}_{/x}}) \ar[r] & \star _{{\mathbf C}_{/y}}}$ dans ${\mathbf C}_{/y}$, d'o\`u en composant avec l'unit\'e $\eta _{\star _{{\mathbf C}_{/y}}}$ de l'adjonction $L_y \dashv U_y$ un morphisme $\xymatrix @C=8mm{f_!(\star _{{\mathbf C}_{/x}}) \ar[r] & U_y (\Omega _y)}$. 
Nous d\'eduisons des adjonctions $f _! \dashv f^*$ et $L_x \dashv U_x$, et de la commutativit\'e du diagramme pr\'ec\'edent les isomorphismes

$$\begin{array}{ccccc}
Hom _{{\mathbf C}_{/y}} ( f_! (\star _{{\mathbf C}_{/x}}) ,U_y (\Omega _y))
& \simeq & Hom  _{{\mathbf C}_{/x}} (\star _{{\mathbf C}_{/x}} , f^* U_y (\Omega _y)) \\
& = & Hom _{{\mathbf C}_{/x}} (\star _{{\mathbf C}_{/x}} , U_x f^*_{[Ab]} (\Omega _y)) & \simeq & Hom _{Ab({\mathbf C}_{/x})} (\Omega _x , f^*_{[Ab]} (\Omega _y)) 
\end{array}
$$
d'o\`u un morphisme naturel 
$\xymatrix@C=8mm{ [{\delta} _f] : \Omega _x \ar[r] & f_{[Ab]}^* (\Omega _y)}$ dans $Ab({\mathbf C}_{/x})$. 

Si le foncteur $f^* _{[Ab]}$ admet un adjoint \`a gauche 
$\xymatrix{{f_!}_{[Ab]}: Ab({\mathbf C}_{/x}) \ar[r] & Ab({\mathbf C}_{/y})}$ 
nous en d\'eduisons un morphisme 
$\xymatrix@C=10mm{[\tilde{\delta _f}] : {f_!}_{[Ab]} (\Omega _x)  \ar[r] & \Omega _y}$ dans $Ab({\mathbf C}_{/y})$, 
qui peut \^etre construit directement \`a partir du diagramme commutatif 
$$\xymatrix @C=10mm{ 
Ab(\mathbf{C}_{/x}) \ar[r] ^{{f_!}_{[Ab]}} & Ab(\mathbf{C}_{/y})  \\ 
 \mathbf{C}_{/x} \ar[u] _{L_x} \ar[r] ^{f_!} & \mathbf{C}_{/y}  \ar[u]_{L_y}
}$$ 

\begin{remark}\label{rmq:existence-adjoint} 
Si la cat\'egorie $\mathbf{C}$ est une cat\'egorie localement pr\'esentable suffisamment r\'eguli\`ere, pour tout objet $x$ dans $\mathbf{C}$ l'adjonction $\xymatrix @C=12mm{
L_x : {\mathbf C}_{/x} \ar@<1ex>[r] & Ab({\mathbf C}_{/x}):U_x \ar@<1ex>[l]-_{\bot}}$
est monadique et la cat\'egorie $Ab(\mathbf{C}_{/x})$ est encore localement pr\'esentable. 

Nous d\'eduisons alors du th\'eor\`eme {\bf 4.5.6.} de \cite{Bo} que le foncteur $f^* _{[Ab]}$ admet toujours un adjoint \`a gauche 
$\xymatrix{{f_!}_{[Ab]}: Ab({\mathbf C}_{/x}) \ar[r] & Ab({\mathbf C}_{/y})}$.  
\end{remark}

\begin{remark}\label{rmq:adjoint} 
Le morphisme $[\tilde{\delta _f}]$ est \'egal au morphisme $L_y(\vartheta _f)$ o\`u $\vartheta _f$ est le morphisme canonique $\xymatrix @C=8mm{\vartheta _f : f_!(\star _{{\mathbf C}_{/x}}) \ar[r] & \star _{{\mathbf C}_{/y}}}$ dans ${\mathbf C}_{/y}$. 
\end{remark}

Le morphisme $[{\delta} _f]$ induit un morphisme $\xymatrix@C=8mm{ U_x( \Omega _x ) \ar[r] & U_x f_{[Ab]}^* (\Omega _y) = f^* U_y (\Omega _y)}$ dans ${\mathbf C}_{/x}$, c'est-\`a-dire un morphisme $\xymatrix@C=8mm{\delta _f: \omega _x \ar[r] & x \times _y \omega _y}$
au dessus de $x$ compatible avec les structures d'objets en groupe ab\'elien de $\omega _x$ et de $x \times _y \omega _y$. 
Alors le morphisme $\xymatrix@C=8mm{\psi _f: x \ar[r] & x \times _y \omega _y}$ dans $Hom _{{\mathbf C}_{/x}} (\star _{{\mathbf C}_{/x}} , U_x f^*_{[Ab]} (\Omega _y))$ correspondant par l'adjonction $L_x \dashv U_x$ au morphisme $[ \delta _f]$ est d\'efini par 
$\psi _f = \delta _f \circ \eta _x$. 

Par construction ce morphisme $\psi _f$ est d\'efini par le diagramme commutatif suivant 
 
$$
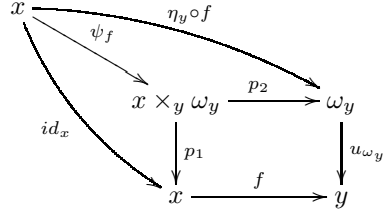
\xymatrix @R=8mm @C=12mm{ x \ar[rd] ^-{\psi _f} \ar@/^3mm/[rrd] ^-{\eta _y \circ f} \ar@/_3mm/[rdd] _-{id _x} \\ 
 & x \times _y \omega _y \ar[r]^{p_2} \ar[d]^-{p_1} & \omega _y \ar[d]^-{u_{\omega _y}} \\  
& x \ar[r]^-{f} & y 
}$$ 
et comme nous avons l'\'egalit\'e $id_y = u_{\omega _y} \circ \eta _y$ nous en d\'eduisons le diagramme suivant form\'e de carr\'es cart\'esiens: 

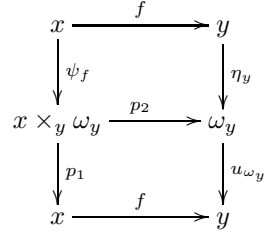
$$\xymatrix @R=8mm @C=12mm{ x \ar[d] ^-{\psi _f} \ar[r]^-f & y \ar[d] ^-{\eta _y}  \\ 
  x \times _y \omega _y \ar[r]^{p_2} \ar[d]^-{p_1} & \omega _y \ar[d]^-{u_{\omega _y}} \\  
x \ar[r]^-{f} & y 
}$$ 
et nous pouvons ainsi \'ecrire $\delta _f \circ \eta _x = \psi _f = f^* (\eta _y)$. 

\vskip .5cm 

Pour tout $[b]$ dans $Ab({\mathbf C}_{/y})$ le morphisme $\xymatrix@C=8mm{ [\delta _f] : \Omega _x \ar[r] & f_{[Ab]}^* (\Omega _y)}$ d\'efinit un morphisme de groupes ab\'eliens 
$$\xymatrix @R=1mm @C=8mm{ 
[\Delta _f] : Der _{\mathbf C}(y,[b]) = Hom _{Ab({\mathbf C}_{/y})} (  \Omega _y , [b])  \ar[r] & Der _{\mathbf C}(x, f^*_{[Ab]} ([b])) = Hom _{Ab({\mathbf C}_{/x})} (  \Omega _x , f^*_{[Ab]} ([b])) \\ 
 {\phantom{\Delta _f : Der _{\mathbf C}(y,[b]) = Hom _{Ab({\mathbf C}_{/y})} } \quad [\delta ]  \phantom{a}} \ar@{|->}[r] & {\phantom{a} f^*_{[Ab]} ([\delta ]) \circ [\delta _f] \quad \phantom{= Hom _{Ab({\mathbf C}_{/x})} (  \Omega _x , f^*_{[Ab]} ([b])) }} 
}$$

Le morphisme dans $Hom _{{\mathbf C}_{/y}} ( \star _{{\mathbf C} _{/y}}, b)$ associ\'e \`a $[\delta ]$ par l'adjonction $L_y \dashv U_y$ est \'egal \`a $\beta = \delta \circ \eta _y$, o\`u $\delta = U_y([\delta ])$, et de m\^eme le morphisme dans $Hom _{{\mathbf C}_{/x}} ( \star _{{\mathbf C} _{/x}}, f^*(b))$ associ\'e \`a $f^*_{[Ab]} ([\delta ]) \circ [\delta _f]$ par l'adjonction $L_x \dashv U_x$ est \'egal \`a $f^*(\delta ) \circ \delta _f \circ \eta _x$. 
Nous d\'eduisons alors de ce qui pr\'ec\`ede que l'application associ\'ee \`a $[\Delta _f]$ par les adjonctions $L_y \dashv U_y$ et $L_x \dashv U_x$ est l'application induite par le foncteur $f^*$: 

$$\xymatrix @R=1mm @C=10mm{ 
\Delta _f : Hom _{{\mathbf C}_{/y}} ( \star _{{\mathbf C} _{/y}}, b) \ar[r] & Hom _{{\mathbf C}_{/x}} ( \star _{{\mathbf C} _{/x}}, f^* (b)) \\ 
{ \phantom{ \Delta _f : Hom _{{\mathbf C}_{/y}} ( \star _{{\mathbf C} _{/y}} } \beta \phantom{a}} \ar@{|->}[r] & {\phantom{a} f^* ( \beta ) \phantom{( \star _{{\mathbf C} _{/x}}, f^* (b))}}
}$$

\vskip .5cm 

\begin{proposition} \label{prop:noyau}
L'\'el\'ement $[\delta ]$ de $Der _{\mathbf C}(y,[b])$ appartient au noyau de $[\Delta _f]$ si et seulement si nous avons l'\'egalit\'e $e_b \circ f = \beta \circ f$, o\`u $\beta$ est le morphisme de $Hom _{{\mathbf C}_{/y}}(\star _{{\mathbf C}_{/y}},b)$ associ\'e \`a $[\delta ]$ et o\`u $e_b$ est l'\'el\'ement neutre de la structure de groupe ab\'elien de $[b]$. 
\end{proposition} 

\begin{preuve} 
Soit $[b]$ un objet de $Ab({\mathbf C}_{/y})$ correspondant \`a un objet $\xymatrix @C=8mm{u_b = b \ar[r] & y}$ de ${\mathbf C}_{/y}$, avec des morphismes $\xymatrix @C=8mm{m_b: b \times _y b \ar[r] & b}$, $\xymatrix @C=8mm{i_b : b \ar[r] & b}$ et $\xymatrix @C=8mm{e_b: y \ar[r] & b}$. 
Alors l'image $[a] = f^*_{[Ab]}([b])$ correspond \`a l'objet $\xymatrix @C=8mm{p_1 : a= f^*(b)=x\times _y b \ar[r] & x}$ de ${\mathbf C}_{/x}$ et aux morphismes $\xymatrix @C=8mm{m_a: a \times _x a \ar[r] & a}$, $\xymatrix @C=8mm{i_a : a \ar[r] & a}$ et $\xymatrix @C=8mm{e_a: x \ar[r] & a}$ obtenus en appliquant le foncteur exact \`a gauche $f^*$. En particulier nous pouvons d\'ecrire le morphisme $e_a$ gr\^ace au diagramme commutatif suivant    

$$\xymatrix @R=6mm @C=12mm{ 
x \ar@/^3mm/[rrd]^-{e_b \circ f} \ar@/_3mm/[rdd]_-{id _x} \ar[rd]^{e_a} \\ 
& x \times _y b \ar[r]^{p_2} \ar[d]^{p_1} & b \ar[d]^{u_b} \\ 
& x \ar[r]^f & y
}$$

D'apr\`es la proposition \ref{prop:morphisme-nul} et la remarque pr\'ec\'edente l'image de $[\delta ]$ est nulle si et seulement si nous avons $f^*(\beta ) = e_a$, c'est-\`a-dire si et seulement si nous avons $e_b \circ f = \beta \circ f$. 

\hfill \tf 
\end{preuve} 

\vskip .5cm 

Le module des diff\'erentielles $\Omega _y$ est le module des diff\'erentielles de $y$ relativement \`a l'objet initial de $\mathbf C$, et nous d\'efinissons le module des diff\'erentielles relatives $\Omega _f$ comme le module des diff\'erentielles de $f$ dans la cat\'egorie des modules de Beck de $f$ vu comme un objet de la cat\'egorie ${\mathbf D} = {\mathbf C}_{x/}$.  
Plus pr\'ecis\'ement pour un morphisme $\xymatrix @C=8mm{f :x \ar[r] & y}$ dans $\mathbf C$, nous consid\'erons $f$ comme un objet de la cat\'egorie ${\mathbf D} = {\mathbf C}_{x/}$, comme pr\'ec\'edemment nous avons le couple de foncteurs adjoints 
$$\xymatrix @C=12mm{
L_f : {\mathbf D}_{/f} \ar@<1ex>[r] & Ab({\mathbf D}_{/f}):U_f \ar@<1ex>[l]-_{\bot}}$$ 
et nous avons la d\'efinition suivante. 

\begin{definition}\label{def:module-des-differentielles-relatives} 
Le \emph{module des diff\'erentielles relatives} $\Omega _f$ est le module de Beck au-dessus de $f$ dans $Ab({\mathbf D}_{/f})$ obtenu comme image par $L_f$ de l'objet final $\star _{{\mathbf D}_{/f}} = (\xymatrix @C=8mm {f \ar[r] ^-{id_f} & f})$ de $\mathbf{D}_{/f}$, 
$\Omega _f = L_f ( \star _{{\mathbf D}_{/F}})$.  
\end{definition}

Nous consid\'erons alors la paire de foncteurs adjoints 
$\xymatrix{ F= [x]_* : {\mathbf C} \ar@<1ex>[r] & {\mathbf C}_{x/} = \mathbf{D}: [x]^! =G \ar[l] }$ 
o\`u l'image de $\xymatrix @C=8mm{(f:x \ar[r] & y)}$ par le foncteur $G$ est $y$, et nous en d\'eduisons la paire de foncteurs adjoints 
$$\xymatrix{ F_y: {\mathbf C} _{/y} \ar@<1ex>[r] &  \mathbf{D} _{/f}: G_f \ar[l] } \ . $$

Les objets de la cat\'egorie ${\mathbf D} _{/f} = ({\mathbf C}_{x/}) _{/f}$ sont les diagrammes commutatifs dans $\mathbf C$: 
$$\xymatrix @R=2mm{ 
& z \ar[dd]^{v} \\ 
x \ar[ur]^{u} \ar[dr]_{f} \\ 
& y
}$$
que nous notons $(u,z,v)$.  
Un morphisme de $(u_1,z_1,v_1)$ dans $(u_2,z_2,v_2)$ est un morphisme $\xymatrix @C=8mm{g:z_1 \ar[r] & z_2}$ 
v\'erifiant $u_2 = g \circ u_1$ et $v_1 = v_2 \circ g$, 
et l'objet final $\star _{{\mathbf D}_{/f}}$ de ${\mathbf D} _{/f}$ est le diagramme $(f,y,id_y)$.  
La donn\'ee d'un morphisme de $\star _{{\mathbf D}_{/f}}$ dans l'objet $(u,z,v)$ est la donn\'ee d'un morphisme $\xymatrix @C=8mm{e:y \ar[r] & z}$ dans $\mathbf C$ v\'erifiant $v\circ e = id_y$ et $e\circ f =u$.

De m\^eme les objets de la cat\'egorie ${\mathbf C} _{/y}$ sont les morphismes $\xymatrix @C=8mm{w: t \ar[r] & y}$ dans $\mathbf C$
que nous notons $(w,t)$.  
Un morphisme de $(w_1,t_1)$ dans $(w_2,t_2)$ est un morphisme $\xymatrix @C=8mm{h: t_1 \ar[r] & t_2}$ 
v\'erifiant $w_1 = w_2 \circ h$, et l'objet final $\star _{{\mathbf C}_{/y}}$ de ${\mathbf C}_{/y}$ est le morphisme $(id_y,y)$. 
La donn\'ee d'un morphisme de $\star _{{\mathbf C}_{/y}}$ dans l'objet $(w,t)$ est la donn\'ee d'un morphisme $\xymatrix @C=8mm{e:y \ar[r] &} t$ dans $\mathbf C$ v\'erifiant $w\circ e = id_y$. 

\vskip .2cm 

Le foncteur $\xymatrix @C=8mm{G_f: {\mathbf D}_{/f} \ar[r] & {\mathbf C}_{/y}}$ induit un foncteur ${{G_f}_{[Ab]}}$ entre les sous-cat\'egories des objets en groupe ab\'eliens, c'est-\`a-dire nous avons le diagramme commutatif suivant: 
$$\xymatrix @R=10mm @C=15mm{ 
Ab(\mathbf{D}_{/f}) \ar[r] ^-{{G_f}_{[Ab]}} \ar[d] ^{U_f} & Ab(\mathbf{C}_{/y}) \ar[d] ^{U_{y}}  \\ 
 \mathbf{D}_{/f} \ar[r] ^-{G_f} & \mathbf{C}_{/y}  
}$$

\vskip .2cm 

\begin{proposition} \label{prop:fidelite-et-equivalence}
\begin{enumerate} 
\item Le foncteur $\xymatrix @C=8mm{{G_f}: {\mathbf D}_{/f} \ar[r] & {\mathbf C}_{/y}}$ est fid\`ele. 
\item Le foncteur $\xymatrix @C=8mm{{{G_f}_{[Ab]}}: Ab({\mathbf D}_{/f}) \ar[r] & Ab({\mathbf C}_{/y})}$ est une \'equivalence de cat\'egories. 
\end{enumerate}
\end{proposition} 

\begin{preuve}   
\begin{enumerate} 
\item 
Il suffit de remarquer que le foncteur $\xymatrix @C=8mm{G_f: {\mathbf D}_{/f} \ar[r] & {\mathbf C}_{/y}}$ envoie $(u,z,v)$ sur $(v,z)$, et un morphisme $\xymatrix @C=8mm{g:z_1 \ar[r] & z_2}$ sur lui-m\^eme. 

\item
Nous allons construire le foncteur ${F_y}_{[Ab]}$ inverse du foncteur ${{G_f}_{[Ab]}}$. 

Soit $(w,t)$ un \'el\'ement de $Ab({\mathbf C}_{/y})$, il existe un morphisme $\xymatrix @C=8mm{e: {\star}_{{\mathbf C}_{/y}} \ar[r] & ( w,t)}$, donc un morphisme $\xymatrix @C=8mm{e: y \ar[r] & t}$ dans $\mathbf C$ v\'erifiant $w\circ e = id_y$. 
Nous d\'efinissons alors ${F_y}_{[Ab]}(w,t)  = (u,t,w)$ avec $u= e \circ f$. 
Il suffit de v\'erifier alors que $(u,t,w)$ est un objet en groupe ab\'elien dans ${\mathbf D}_{/f}$, et que nous avons ${{G_f}_{[Ab]}} \circ {{F_y}_{[Ab]}} = id _{Ab({\mathbf C}_{/y})}$ et ${{F_y}_{[Ab]}} \circ {{G_f}_{[Ab]}} = id _{Ab({\mathbf D}_{/f})}$.  
\end{enumerate}
\hfill$\Box$ 
\end{preuve}    

\vskip .2cm 

\begin{remark}\label{rmq:objet-abelien-df}
Nous pouvons d\'ecrire tout objet en groupe dans ${\mathbf D}_{/f}$ par un diagramme $\xymatrix @C=8mm{x \ar[r]^-{v_z} & z \ar[r]^-{u_z} & y \ar@/^/[l]^-{e_z}}$, avec $u_z \circ e_z = id _y$, $u_z \circ v_z =f$  et $e_z \circ f =v_z$, o\`u le morphisme $e_z$ est d\'efini par l'\'el\'ement neutre. 
En particulier le module des diff\'erentielles $\Omega _f$ correspond \`a un diagramme 
$\xymatrix @C=8mm{x \ar[r]^-{v_{\omega _f}} & \omega _f \ar[r]^-{u_{\omega _f}} & y \ar@/^/[l]^-{e_{\omega _f}}}$. 
\end{remark} 

\vskip .2cm 

\begin{lemma}\label{lem:transformation-naturelle}  
Il existe une transformation naturelle $[\Gamma ] : L_yG_f \Rightarrow {G_f}_{[Ab]}L_f$ telle que pour tout $a$ dans ${\mathbf D}_{/f}$ le morphisme $[\Gamma (a)] \in Hom _{Ab({\mathbf C}_{/x})} (L_yG_f (a) , {G_f}_{[Ab]}L_f (a))$ est un \'epimorphisme.

\end{lemma}  

\begin{preuve} 
La transformation naturelle $[\Gamma ]$ est d\'efinie par $[\Gamma ] =  \epsilon _y  {G_f}_{[Ab]} L_f \circ  L_y G_f \eta _f$ o\`u $\eta _f$ est l'unit\'e de l'adjonction $L_f \dashv U_f$, $\epsilon _y$ la counit\'e de l'adjonction $L_y \dashv U_y$ et gr\^ace \`a l'\'egalit\'e $G_f U_f = U_y {G_f}_{[Ab]}$.  

Pour tout $a$ dans ${\mathbf D}_{/f}$ et tout $[b]$ nous d\'eduisons de l'adjonction $L_f \dashv U_f$ et de la pleine fid\'elit\'e du foncteur ${G_f}_{[Ab]}$ les isomorphismes 
$$Hom _{{\mathbf D}_{/f}}(a, U_f([b])) \simeq  Hom _{Ab({\mathbf D}_{/f})}(L_f(a), [b]) \simeq Hom _{Ab({\mathbf C}_{/y})}({G_f}_{[Ab]}L_f(a), {{G_f}_{[Ab]}}([b])) \ ,$$
et de la commutativit\'e du diagramme pr\'ec\'edent et de l'adjonction $L_y \dashv U_y$ les isomorphismes 
$$ Hom _{{\mathbf C}_{/y}}(G_f(a), G_fU_f([b])) \simeq Hom _{{\mathbf C}_{/y}}(G_f(a), U_y {G_f}_{[Ab]}([b]))  \simeq Hom _{Ab({\mathbf C}_{/y})}(L_yG_f(a), {{G_f}_{[Ab]}}([b])) \ .$$
Le foncteur $G_f$ induit une application de $Hom _{{\mathbf D}_{/f}}(a, U_f([b]))$ dans $Hom _{{\mathbf C}_{/y}}(G_f(a), G_fU_f([b]))$, d'o\`u une application 
$$\gamma (a) : \xymatrix{ Hom _{Ab({\mathbf C}_{/y})}({G_f}_{[Ab]}L_f(a), {{G_f}_{[Ab]}}([b])) \ar[r] & Hom _{Ab({\mathbf C}_{/y})}(L_yG_f(a), {{G_f}_{[Ab]}}([b])) \ .}$$ 

Alors l'application $\gamma (a)$ correspond \`a la composition par le morphisme $[\Gamma (a)]$. 
De la fid\'elit\'e du foncteur $G_f$ et de l'essentielle surjectivit\'e du foncteur ${G_f}_{[Ab]}$ nous d\'eduisons que pour tout $[c]$ dans $Ab({\mathbf C}_{/y})$  l'application  
$$ \xymatrix{ Hom _{Ab({\mathbf C}_{/y})}({G_f}_{[Ab]}L_f(a), [c]) \ar[r] & Hom_{Ab({\mathbf C}_{/y})}(L_yG_f(a), [c])}$$ 
induite par le morphisme $[\Gamma (a)]$ est injective, par cons\'equent $[\Gamma (a)]$ est un \'epimorphisme. 

\hfill \tf 
\end{preuve} 

\vskip .5cm 

Pour $a$ \'egal \`a l'objet final $\star _{{\mathbf D}_{/f}}$, nous trouvons un \'epimorphisme naturel:  
$$[\gamma _f] : \xymatrix{\Omega _y \ar[r] & {G_f}_{[Ab]} ( \Omega _f)} \ ,$$ 
et nous d\'eduisons que pour tout $[b]$ dans $Ab({\mathbf D}_{/f})$ nous avons un morphisme injectif entre les ensembles des d\'erivations $\xymatrix @C=8mm{ [\Gamma _f] : Der _{\mathbf D}(f, [b]) \ar@<-2pt>@{^{(}->}[r] & Der_{\mathbf C} (y, {G_f}_{[Ab]}([b]))}$.  

 \begin{remark} 
 Le morphisme $[\gamma _f]$ se factorise $\xymatrix@C=8mm{\Omega _y \ar[r] & L_yU_y {G_f}_{[Ab]} (\Omega _f) \ar[r] & {G_f}_{[Ab]} ( \Omega _f)}$.  
 \end{remark} 
 
\vskip .5cm 

\begin{definition} 
Soient $\xymatrix @C=8mm{f : X \ar[r] &Y}$ et $\xymatrix @C=8mm{g : Y \ar[r] &Z}$ deux morphismes dans une cat\'egorie additive $\mathbf A$, alors on dit que la suite 
$$\xymatrix @C=8mm{ X \ar[r] ^R & Y \ar[r] ^S & Z \ar[r] & 0}$$
est \emph{exacte} dans $\mathbf A$ si pour tout $U$ dans $\mathbf A$ nous avons une suite excate de groupes ab\'eliens 
$$\xymatrix @C=8mm{ 
0 \ar[r]  & Hom_{\mathbf A} (Z,U)  \ar[r] ^{S^*} & Hom_{\mathbf A} (Y,U) \ar[r] ^{R^*} & Hom_{\mathbf A} (X,U) \ .
}$$
\end{definition} 

\vskip .5cm 

\begin{theorem} \label{th:premiere-suite-exacte}
Soit $\mathbf C$ une cat\'egorie localement pr\'esentable, et soit $\xymatrix @C=8mm{f:x \ar[r] & y}$ dans $\mathbf C$ tel que le foncteur $f^* _{[Ab]}$ admette un adjoint \`a gauche 
$\xymatrix{{f_!}_{[Ab]}: Ab({\mathbf C}_{/x}) \ar[r] & Ab({\mathbf C}_{/y})}$, alors nous avons une suite exacte dans $Ab({\mathbf C} _{/y})$: 
$$\xymatrix @C=4mm{{f_!}_{[Ab]} ( \Omega _x) \ar@<-1pt>[rrr]^-{[\tilde{\delta _f}]} &&& \Omega _y \ar@<-1pt>[rrr] ^-{[\gamma _f]  } &&& {G_f}_{[Ab]}  (\Omega _f) \ar[rr] && 0}$$ 
\end{theorem} 

\begin{preuve}   
Par d\'efinition il faut montrer que pour tout $[c]$ dans $Ab({\mathbf C} _{/y})$ la suite de morphismes induite 
$$\xymatrix @C=5mm {0 \ar[r] & Hom_{Ab({\mathbf C}_{/y})} ({G_f}_{[Ab]}  ( \Omega _f), [c]) \ar@<-1pt>[rr]^-{[\gamma _f]^*}  && Hom_{Ab({\mathbf C}_{/y})} (\Omega _y , [c]) \ar@<-1pt>[rr]^-{[\tilde\delta _f]^*} && Hom_{Ab({\mathbf C} _{/y})} ( {f_!}_{[Ab]}(\Omega _x), [c])}$$
est une suite exacte de groupes ab\'eliens. 
De ce qui préc\`ede nous en d\'eduisons que c'est \'equivalent \`a monter que pour tout $[b]$ dans $Ab({\mathbf D} _{/f})$ la suite de morphismes 
$$\xymatrix @C=5mm{ 0 \ar[r] & Hom_{Ab({\mathbf D}_{/f})} (\Omega _f, [b]) \ar@<-1pt>[rr]^-{[\Gamma _f]} && Hom_{Ab({\mathbf C}_{/y})} (\Omega _y , {G_f}_{[Ab]}([b])) \ar@<-1pt>[rr]^-{[\Delta _f]} && Hom_{Ab({\mathbf C} _{/x})} (\Omega _x, f_{[Ab]}^* {G_f}_{[Ab]}([b]))}$$ 
est exacte. 

Comme nous avons d\'ej\`a vu que le morphisme $[\Gamma _f]$ est injectif il suffit de monter que nous avons l'\'egalit\'e: 
$$Ker( [\Delta _f] ) = Im ( [\Gamma _f]) \ .$$ 
Nous pouvons \'ecrire le diagramme commutatif suivant 
$$\xymatrix @C=10mm @R=7mm{Der _{\mathbf D}(f, [b]) \ar@<-1pt>[r]^-{[\Gamma _f]} \ar@{=}[d] & Der_{\mathbf C} (y, {G_f}_{[Ab]}([b])) \ar@<-1pt>[r]^-{[\Delta _f]}  \ar@{=}[d]  & Der_{\mathbf C} (x, f_{[Ab]}^*{G_f}_{[Ab]}([b]))  \ar@{=}[d] \\ 
Hom_{Ab({\mathbf D}_{/f})} (\Omega _f, [b]) \ar@<-1pt>[r]^-{[\Gamma _f]} \ar[d]^(.45)*[@]{\simeq} & Hom_{Ab({\mathbf C}_{/y})} (\Omega _y , {G_f}_{[Ab]}([b])) \ar@<-1pt>[r]^-{[\Delta _f]} \ar[d]^(.45)*[@]{\simeq} & Hom_{Ab({\mathbf C} _{/x})} (\Omega _x, f_{[Ab]}^* {G_f}_{[Ab]}([b])) \ar[d]^(.45)*[@]{\simeq} \\ 
Hom_{{\mathbf D}_{/f}} (\star_{{\mathbf D}_{/f}}  , U_f([b])) \ar@<-1pt>[r]^-{\Gamma _f}& Hom_{{\mathbf C}_{/y}} (\star _{{\mathbf C}_{/y}} , U_y{G_f}_{[Ab]}([b])) \ar@<-1pt>[r] ^-{\Delta _f} & Hom_{{\mathbf C} _{/x}} (\star _{{\mathbf C} _{/x}}  , U_x f_{[Ab]}^* {G_f}_{[Ab]}([b])) 
}$$

D'apr\`es la remarque \ref{rmq:objet-abelien-df} un objet $[b]$ de $Ab({\mathbf D}_{/f})$ correspond \`a un diagramme $\xymatrix @C=8mm{x \ar[r]^-{v_b} & b \ar[r]^-{u_b} & y }$, avec $\xymatrix @C=8mm{e_b  : y \ar[r] & b}$ v\'erifiant $e_b \circ f =v_b$. 
Un morphisme $[\alpha] \in Hom_{{\mathbf D}_{/f}} (\star_{{\mathbf D}_{/f}}  , U_f([b]))$ correspond \`a un morphisme $\xymatrix @C=8mm{\alpha : y \ar[r] &b}$ tel que $\alpha \circ f = v_b$ et $u_b \circ \alpha = id _y$, 
et un morphisme $[\beta] \in Hom_{{\mathbf C}_{/y}} (\star _{{\mathbf C}_{/y}} , U_y{G_f}_{[Ab]}([b]))$ correspond \`a un morphisme $\xymatrix @C=8mm{\beta : y \ar[r] &b}$ tel que $u_b \circ \beta = id _y$, et l'application $\Gamma _f$ envoie le morphisme $[\alpha]$ sur $[\beta]$ avec $\alpha = \beta$. 

En particulier un morphisme $[\beta ]$ de $Hom_{{\mathbf C}_{/y}} (\star _{{\mathbf C}_{/y}} , U_y{G_f}_{[Ab]}([b]))$ appartient \`a l'image de l'application $\Gamma _f$ si et seulement si nous avons $\beta \circ f= v_b$, c'est-\`a-dire si et seulement si nous avons $\beta \circ f= e_b \circ f$. Le r\'esultat est alors une cons\'equence de la proposition \ref{prop:noyau}.  

\hfill$\Box$ 
\end{preuve}  

\vskip .5cm 

Sous certaines conditions le module des diff\'erentielles relatif $\Omega _f$ est nul; plus pr\'ecis\'ement nous avons le r\'esultat suivant. 

\begin{theorem} \label{th:deuxieme-suite-exacte} 
Soit $\mathbf C$ une cat\'egorie localement pr\'esentable, et soit $\xymatrix @C=8mm{f:x \ar[r] & y}$ dans $\mathbf C$ tel que le foncteur $f^* _{[Ab]}$ admette un adjoint \`a gauche 
$\xymatrix{{f_!}_{[Ab]}: Ab({\mathbf C}_{/x}) \ar[r] & Ab({\mathbf C}_{/y})}$, alors si $f$ est un \'epimorphisme dans $\mathbf C$, l'application 
$\xymatrix @C=8mm{ [\tilde{\delta _f}] :{f_!}_{[Ab]} ( \Omega _x) \ar[r] & \Omega _y}$ 
est un \'epimorphisme dans $Ab({\mathbf C}_{/y})$.
\end{theorem} 

\begin{preuve}   
Si $\xymatrix @C=8mm{f:x \ar[r] & y}$ est un \'epimorphisme dans $\mathbf C$, le morphisme induit $\xymatrix @C=8mm{\vartheta _f : f_!(\star _{{\mathbf C}_{/x}}) \ar[r] & \star _{{\mathbf C}_{/y}}}$, correspondant \`a $\xymatrix @C=8mm{\vartheta _f : ( x \ar[r] ^f & y ) \ar[r] & (y \ar[r] ^{id_y} & y)}$ est un \'epimorphisme dans ${\mathbf C}_{/y}$. 
Comme le foncteur $L_y$ est un adjoint \`a gauche, il en est de m\^eme du morphisme $L_y(\vartheta _f)$, qui est \'egal au morphisme $[\tilde{\delta _f}]$ d'apr\`es la remarque \ref{rmq:adjoint}. 

\hfill \tf 
\end{preuve} 

\vskip .5cm 

Nous consid\'erons deux objets $x$ et $z$ de $\mathbf C$, $y= z \coprod x$, et les morphismes $\xymatrix @C=8mm{f=q_2:x \ar[r] & y}$ et $\xymatrix @C=8mm{g=q_1:z \ar[r] & y}$. 
Comme pr\'ec\'edemment nous avons la paire de foncteurs adjoints 
$$\xymatrix{ F^{(1)} = [z]_* : {\mathbf C} \ar@<1ex>[r] & {\mathbf C}_{z/} = \mathbf{E}: [z]^! = G^{(1)} \ar[l] } \ ,$$ 
l'image de $\xymatrix @C=8mm{(g:z \ar[r] & y)}$ par le foncteur $G^{(1)}$ est $y$, nous en d\'eduisons la paire de foncteurs adjoints 
$$\xymatrix{ F^{(1)}_y: {\mathbf C} _{/y} \ar@<1ex>[r] &  \mathbf{E} _{/g}: G^{(1)}_g \ar[l] }\ . $$  

Comme pr\'ec\'edemment nous supposons que le foncteur $\xymatrix @C=8mm{f^* _{[Ab]} : Ab({\mathbf C}_{/y}) \ar[r] & Ab({\mathbf C}_{/x})}$ admet un adjoint \`a gauche ${f_!}_{[Ab]}$.

\begin{theorem} \label{th:troisieme-suite-exacte}  
Il existe dans $Ab({\mathbf C}_{/y})$ un isomorphisme naturel
 $$ {f_!}_{[Ab]} (\Omega _x) \simeq {G_g}_{[Ab]}^{(1)} (\Omega _g) \ . $$
\end{theorem} 

\begin{preuve} 
Le foncteur $\xymatrix @C=8mm{{G_g}_{[Ab]}^{(1)}: Ab({\mathbf E}_{/g}) \ar[r] & Ab({\mathbf C}_{/y})}$ est une \'equivalence de cat\'egories, 
le foncteur inverse ${F_y}_{[Ab]}^{(1)}$ est l'adjoint \`a gauche et nous avons le diagramme commutatif suivant

$$\xymatrix @C=12mm @R=10mm{ 
Ab(\mathbf{C}_{/x}) \ar[r] ^{{f_!} _{[Ab]}} & Ab(\mathbf{C}_{/y}) \ar[r] ^{{F_y} _{[Ab]}^{(1)}} & Ab(\mathbf{E}_{/g})  \\ 
\mathbf{C}_{/x} \ar[u] _{L_x} \ar[r]  ^{f_!} & \mathbf{C}_{/y} \ar[u] _{L_y} \ar[r]  ^{F^{(1)}_y} & \mathbf{E}_{/g} \ar[u]_{L_g}
}$$

Par d\'efinition du foncteur $f_!$ nous avons 
$f_! ( {\star}_{{\mathbf C}_{/x}}) = \xymatrix @C=8mm{(f: x \ar[r] & y)}$ 
et d'apr\`es la remarque \ref{rmq:foncteurs-induits} nous avons   $F^{(1)}_y \xymatrix @C=8mm{(f: x \ar[r] & y)} = {\star}_{{\mathbf E}_{/g}}$, 
nous d\'eduisons alors de la commutativit\'e du diagramme pr\'ec\'edent l'\'egalit\'e 
$${F_y}_{[Ab]}^{(1)}{f_!}_{[Ab]}(\Omega _x) = \Omega _g\ ,$$ 
d'o\`u l'isomorphisme cherch\'e. 

\hfill \tf 
\end{preuve} 

\vskip .5cm


\section{Exemples} 


\subsection{Cat\'egorie des ensembles}  

Nous prenons comme cat\'egorie $\mathbf C$ la cat\'egorie $\mathbf{Ens}$ des ensembles.   
La cat\'egorie des objets en groupe ab\'elien dans $\mathbf{Ens}$ est la cat\'egorie $\mathbf{Ab}$ des groupes ab\'eliens. 

Pour tout $X$ dans $\mathbf{Ens}$  la cat\'egorie $Ab(\mathbf{Ens}_{/X})$ est la cat\'egorie des applications  $\xymatrix @C=8mm{u : F \ar[r] & X}$ telles que chaque fibre $F_x = u^{-1}(x)$ pour $x \in X$ est un groupe ab\'elien, et les morphismes sont d\'efinis par 
$$Hom _{Ab(\mathbf{Ens}_{/X})} ( (\xymatrix @C=8mm{F^{(1)} \ar[r] ^-{u^{(1)}} & X}) , (\xymatrix @C=8mm{F^{(2)} \ar[r] ^-{u^{(2)}} & X} )) = \prod _{x \in X} Hom _{\mathbf{Ab}} ( F^{(1)}_x ,  F^{(2)}_x ) \ . $$ 
Le foncteur $\xymatrix @C=8mm{L_X : \mathbf{Ens}_{/X} \ar[r] & Ab(\mathbf{Ens}_{/X})}$ adjoint \`a gauche du foncteur  $\xymatrix @C=8mm{U_X : Ab(\mathbf{Ens}_{/X}) \ar[r] & \mathbf{Ens}_{/X}}$ est d\'efini par $L_X (\xymatrix @C=8mm{E \ar[r] ^-{l} & X)} = (\xymatrix @C=8mm{F = \mathbb{Z} _X[E] \ar[r] ^-{u} & X)}$ avec $F_x = \mathbb{Z}[E_x]$ le groupe ab\'elien libre engendr\'e par la fibre $E_x$ pour tout $x \in X$.   
Nous en d\'eduisons que le module des diff\'erentielles de l'ensemble $X$ est d\'efini par $\Omega _X = \xymatrix @C=8mm{(X \times \mathbb{Z} \ar[r] ^-{p_1} & X)}$. 

\vskip .5cm 

Soit $\xymatrix @C=8mm{f : X \ar[r] & Y}$ une application dans $\mathbf{Ens}$, alors le foncteur $\xymatrix @C=8mm{f ^* _{[Ab]} :  Ab(\mathbf{Ens}_{/Y}) \ar[r] & Ab(\mathbf{Ens}_{/X})}$ envoie l'objet en groupe ab\'elien $\xymatrix @C=8mm{(G \ar[r] ^-{v} & Y)}$ sur $\xymatrix @C=8mm{(F = X \times _Y G  \ar[r] ^-{p_1} & X)}$ o\`u pour tout $x \in X$ nous avons $F_x = G_{f(x)}$. 

Ce foncteur admet un adjoint \`a gauche $\xymatrix @C=8mm{{f _!} _{[Ab]} :  Ab(\mathbf{Ens}_{/X}) \ar[r] & Ab(\mathbf{Ens}_{/Y})}$ qui envoie l'objet en groupe ab\'elien $\xymatrix @C=8mm{(F \ar[r] ^-{u} & X)}$ sur $\xymatrix @C=8mm{(G  \ar[r] ^-{v} & Y)}$ o\`u pour tout $y \in Y$ nous avons $G_y = \coprod_{f(x) =y} F_x$. 

En effet pour $\xymatrix @C=8mm{(F^{(1)} \ar[r] ^-{u ^{(1)}} & X)}$ dans $Ab(\mathbf{Ens}_{/X})$ et pour $\xymatrix @C=8mm{(G ^{(2)}  \ar[r] ^-{v ^{(2)}} & Y)}$ dans $Ab(\mathbf{Ens}_{/Y})$ nous avons 

$$Hom _{Ab(\mathbf{Ens}_{/X})} ( (\xymatrix @C=8mm{F^{(1)} \ar[r] ^-{u^{(1)}} & X}) , f ^* _{[Ab]} (\xymatrix @C=8mm{G^{(2)} \ar[r] ^-{v^{(2)}} & Y} )) = \prod _{x \in X} Hom _{\mathbf{Ab}} ( F^{(1)}_x ,  G^{(2)}_{f(x)} )  $$ 

et 

$$Hom _{Ab(\mathbf{Ens}_{/Y})} ( {f_!}_{[Ab]} (\xymatrix @C=8mm{F^{(1)} \ar[r] ^-{u^{(1)}} & X}) , (\xymatrix @C=8mm{G^{(2)} \ar[r] ^-{v^{(2)}} & Y} )) = \prod _{y \in Y} Hom _{\mathbf{Ab}} ( \coprod _{f(x) =y} F^{(1)}_x ,  G^{(2)}_y ) \ . $$

Pour une application $\xymatrix @C=8mm{f : X \ar[r] & Y}$ nous avons alors ${f_!} _{[Ab]} ( \Omega _X ) = \xymatrix @C=8mm{(G   \ar[r] ^-{v} & Y)}$ avec pour tout $y \in Y$ la fibre $G_y = \coprod _{f(x) = y} \mathbb{Z}$ et le morphisme $\xymatrix @C=8mm{[\tilde{\delta _f}] : {f_!} _{[Ab]} ( \Omega _X ) \ar[r] & \Omega _Y}$ est le morphisme qui au dessus de $y$  est d\'efini par 
$$\xymatrix @C=8mm{ \displaystyle \delta _y : \coprod _{f(x) = y} \mathbb{Z} \ar[r] & \mathbb{Z}} \ .$$

Si l'ensemble $f^{-1}(y) = \{ x \in X \ | \ f(x)=y \}$ est non vide nous avons $\xymatrix @C=8mm{\delta _y : \coprod _{f(x) = y} \mathbb{Z} \ar[r] & \mathbb{Z}}$ qui est d\'efini par $id _{\mathbb{Z}}$ sur chaque facteur et $\delta _y$ est surjectif, et si l'ensemble $f^{-1}(y)$ est vide nous avons $\xymatrix @C=8mm{ \displaystyle \delta _y : (0) \ar[r] & \mathbb{Z}}$. 
Nous en d\'eduisons que le module des diff\'erentielles relatives est d\'efini par 
$\xymatrix @C=8mm{\Omega _f = ( H \ar[r] & Y )}$ avec 
$H_y = (0)$ si $y$ appartient \`a l'image de l'application $f$ et $H_y = \mathbb{Z}$ sinon. 

Nous d\'eduisons de ce qui pr\'ec\`ede le r\'esultat suivant. 

\begin{proposition}\label{prop:ens} 
Le morphisme $\xymatrix @C=8mm{[\tilde{\delta _f}] : {f_!} _{[Ab]} ( \Omega _X ) \ar[r] & \Omega _Y}$ est un \'epimorphisme (resp. un monomorphisme) si et seulement si l'application $\xymatrix @C=8mm{f : X \ar[r] & Y}$ est surjective (resp. injective). 

En particulier $\xymatrix @C=8mm{[\tilde{\delta _f}] : {f_!} _{[Ab]} ( \Omega _X ) \ar[r] & \Omega _Y}$ est un isomorphisme dans $Ab(\mathbf{Ens}_{/Y})$ si et seulement si $\xymatrix @C=8mm{f : X \ar[r] & Y}$ est une bijection dans $\mathbf{Ens}$. 
\end{proposition} 

\vskip .5cm

\subsection{Cat\'egorie des mono\"\i des commutatifs} 

Nous consid\'erons maintenant la cat\'egorie $\mathbf{Com}$ des mono\"\i des associatifs, commutatifs unitaires. 
Un objet $X$ de $\mathbf{Com}$ est un ensemble muni d'une loi de composition interne $\xymatrix @C=8mm{\star : X \times X \ar[r] & X}$ associative, commutative et muni d'un \'el\'ement neutre $1_X$ pour la loi $\star$. 
Suivant Michael Barr nous avons la description suivante.

Soit $X$ dans $\mathbf{Com}$, un objet de la cat\'egorie $Ab( \mathbf{Com} _{/X})$ est la donn\'ee d'un morphisme de mono\"\i des $\xymatrix @C=8mm{u_A : A \ar[r] & X}$ et de morphismes de mono\"\i des au-dessus de $X$ $\xymatrix @C=8mm{m_A : A \times _X A \ar[r] & A}$, $\xymatrix @C=8mm{i_A : A \ar[r] & A}$ et $\xymatrix @C=8mm{e_A : X \ar[r] & A}$ qui d\'efinissent une structure d'objet en groupe ab\'elien. 

Nous en d\'eduisons que pour tout $x \in X$ ces morphismes induisent une structure de groupe ab\'elien sur la fibre $A_x = u_A^{-1}(x)$ not\'ee $+_x$ dont l'\'el\'ement neutre est $0_{A_x}= e_A(x)$. 
La composition par un \'el\'ement $x$ de $X$ d\'efinit une application de $A_y$ dans $A_{x\star y}$, et comme d'apr\`es la d\'emonstration de la proposition \ref{prop:loi-de-groupe} cette application est un morphisme de groupes ab\'eliens.   
Cela d\'etermine compl\`etement les objets en groupe ab\'elien dans $\mathbf{Com} _{/X}$, plus pr\'ecis\'ement nous avons le r\'esultat suivant. 

\begin{proposition}\cite{Ba2} \label{prop:monoide} 
Soit $X$ un mono\"\i de commutatif, un objet $A$ de $Ab( \mathbf{Com} _{/X})$ est la donn\'ee pour tout $x$ dans $X$ d'un groupe ab\'elien $A_x$, pour tout $x$ et $y$ dans $X$ d'un morphisme de groupes $\xymatrix @C=8mm{ h^{(A)}_x : A_y \ar[r] & A_{x \star y}}$ v\'erifiant $h^{(A)}_x \circ h^{(A)}_y =h^{(A)}_{x \star y}$. 

Nous notons un \'el\'ement de $A$ sous la forme d'une paire $(x,a)$ avec $x=u_A(x,a)$ dans $X$, $a$ dans $A_x$ et la structure de mono\"\i de sur $A$ est d\'efinie par $(x,a) \star (y,b) = \bigl (x \star y , h^{(A)}_x(b) + _{x \star y} h^{(A)}_{y}(a) \bigr )$.  
\end{proposition} 

Pour tout $A$ dans $Ab( \mathbf{Com} _{/X})$ l'espace des d\'erivations $Der_{\mathbf {Com}}(X,A)$ est \'egal par d\'efinition \`a l'ensemble $Hom _{\mathbf{Com} _{/X}}(X,A)$. 
D'apr\`es ce qui pr\'ec\`ede un morphisme $s$ dans $Hom _{\mathbf{Com} _{/X}}(X,A)$ correspond \`a la donn\'ee pour tout $x$ dans $X$ d'un \'el\'ement $s_x$ dans $A_x$ v\'erifiant la propri\'et\'e suivante 
$$s_{x \star y} = \bigr (h^{(A)}_x(s_{y}) + _{x \star y} h^{(A)}_{y}(s_x) \bigr ) \ .$$
En particulier pour $x=y$ \'egal \`a l'\'el\'ement neutre $1_X$ du mono\"\i de $X$, nous avons $h^{(A)}_x = id$ et nous en d\'eduisons $s_{1_X}$ est \'egal \`a l'\'el\'ement neutre $0_{A_{1_X}}$ du groupe ab\'elien $A_{1_X}$. 

\vskip .2cm

Un morphisme $\xymatrix @C=7mm{f : A \ar[r] & B}$ dans $Ab( \mathbf{Com} _{/X})$ induit pour chaque $x$ dans $X$ un morphisme de groupes ab\'eliens $\xymatrix @C=8mm{f_x : A_x \ar[r] & B_x}$. 
Comme $f$ respecte la structure de mono\"\i des nous avons de plus pour tout $a$ dans $A_x$ et $b$ dans $A_{y}$  l'\'egalit\'e 
$f_{x\star y} \bigl ( h^{(A)}_x (b) + _{x \star y}  h^{(A)}_{y} (a) \bigr ) = h^{(B)}_x\bigl ( f_{y} (b) \bigr ) + _{x \star y}  h^{(B)}_{y} \bigl ( f_x (a) \bigr )$, 
ce qui induit pour $b=0_{A_y}$ l'\'egalit\'e $f_{x \star y} ( h^{(A)}_x (a)) =  h^{(B)}_x \bigl ( f_{y}(a) \bigr )$. 
Nous avons alors le r\'esultat suivant. 

\begin{proposition}\cite{Ba2} \label{prop:morphisme-de-monoides} 
Un morphisme $\xymatrix @C=7mm{f : A \ar[r] & B}$ dans $Ab( \mathbf{Com} _{/X})$ est d\'efini par la donn\'ee d'une famille de morphismes de groupes ab\'eliens $\xymatrix @C=8mm{f_x : A_x \ar[r] & B_x}$ v\'erifiant la relation  $f_{x \star y} \circ h^{(A)}_x =  h^{(B)}_x \circ f_{y}$.   
\end{proposition} 

\vskip .2cm 

Soit $\xymatrix @C=8mm{f : X \ar[r] & Y}$ un morphisme dans $\mathbf{Com}$, nous avons alors la paire de foncteurs adjoints 
$$\xymatrix{f_! : {\mathbf{Com}}_{/X} \ar@<1ex>[r] & {\mathbf{Com}}_{/Y} : f^* \ar[l]}$$ 
et le diagramme commutatif 
$$\xymatrix @C=12mm{ 
Ab(\mathbf{Com}_{/Y}) \ar[r] ^-{f_{[Ab]}^*} \ar[d] ^{U_Y} & Ab(\mathbf{Com}_{/X}) \ar[d] ^{U_X} \\ 
 \mathbf{Com}_{/Y} \ar[r] ^-{f^*} & \mathbf{Com}_{/X}
}$$
o\`u le foncteur $f^* _{[Ab]}$ envoie l'objet $B$ de $Ab({\mathbf{Com}}_{/Y})$ d\'efini par la famille de groupes ab\'eliens $\bigl ( B_y \bigr ) _{y\in Y}$ et les morphismes $h^{(B)}_y$ sur l'objet $A= f^* _{[Ab]} (B)$ de $Ab( {\mathbf{Com}}_{/X})$ d\'efini par la famille de groupes ab\'eliens $\bigl ( A_x \bigr ) _{x \in X}$ avec $A_x= B_{f(x)}$ et les morphismes $h^{(A)} _x = h^{(B)} _{f(x)}$. 

\vskip .2cm 

Soit $\Omega _X$ le module des diff\'erentielles de $X$ et soit $\xymatrix @C=8mm{\eta _X : X \ar[r] & \Omega _X}$ l'unit\'e de l'adjonction, alors $\Omega _X$ correspond \`a une famille $\bigl ( ({\Omega _X}) _x \bigr ) _{x \in X}$ de groupes ab\'eliens et \`a des morphismes $\xymatrix @C=6mm{h_x^{(\Omega _X)} :{ (\Omega _X}) _y \ar[r] & ({\Omega _X}) _{x \star y}}$, et le morphisme $\eta _X$ correspond \`a la donn\'ee d'\'el\'ements $\eta _x$ dans $({\Omega _X}) _x$ v\'erifiant 
$$\eta _{x \star y} = \bigr (h^{(\Omega _X)}_x(\eta _{y}) + _{x \star y} h^{(\Omega _X)}_{y}(\eta _x) \bigr ) \ .$$

\vskip .2cm 

Nous consid\'erons le mono\"\i de $X= ({\mathbb N},+)$, un objet $A$ dans $Ab( \mathbf{Com} _{/\N})$ est d\'efini par la donn\'ee de groupes ab\'eliens $A_n$ et de morphismes de groupes $\xymatrix @C=8mm{h^{(A)} : A_n \ar[r] & A_{n+1}}$, pour tout $n$ appartenant \`a ${\mathbb N}$, avec les notations pr\'ec\'edentes nous avons $h^{(A)} _n = \bigl ( h^{(A)} \bigr ) ^n$. 
Et de la m\^eme mani\`ere un morphisme $\xymatrix @C=8mm{f : A \ar[r] & B}$ dans $Ab( \mathbf{Com} _{/\N})$ est d\'efini par une famille de morphismes de groupes ab\'eliens $\xymatrix @C=8mm{f _n: A_n \ar[r] & B_n}$ v\'erifiant $f_{n+1} \circ h^{(A)} =  h^{(B)} \circ f_{n}$.   

\begin{proposition} \label{exemple-monoide}
Le module des diff\'erentielles $\Omega _{\N}$ est l'objet de $Ab( \mathbf{Com} _{/\N})$ d\'efini par 
\begin{enumerate}
\item $({\Omega _{\N}})_0 = (0)$ 

\item $({\Omega _{\N}} )_n = \Z$ et $\xymatrix @C=12mm{h^{(\Omega _{\N})} : ({\Omega _{\N}} )_n  \ar[r] ^{id_{\Z}} & ({\Omega _{\N}} )_{n+1}}$ pour $n \geq 1$. 
\end{enumerate} 
\end{proposition} 

\begin{preuve}

Pour tout $A$ dans $Ab( \mathbf{Com} _{/\N})$ l'ensemble $Hom_{\mathbf{Com} _{/\N}}(\N ,A)$ est isomorphe au groupe ab\'elien $A_1$. 
En effet un \'el\'ement $s_1$ dans $A_1$ d\'etermine un morphisme $\xymatrix @C=8mm{s: \N \ar[r] & A}$ avec $s_0= 0_{A_0}$ et pour tout $n \geq 1$ $s_n = n \bigl ( h^{(A)}\bigr ) ^{n-1} (s_1)$. 

Un morphisme $\xymatrix @C=8mm{f : \Omega _{\N} \ar[r] & A}$ dans $Ab(\mathbf{Com} _{/\N})$ est enti\`erement d\'etermin\'e par $\xymatrix @C=8mm{f_1 : ({\Omega _{\N}})_1 \ar[r] & A_1}$, 
en effet nous avons $f_0=0$ et $f_n = \bigl ( h^{(A)} \bigr ) ^n \circ f_1$ pour tout $n \geq 1$. 
Nous d\'eduisons de l'\'egalit\'e $({\Omega _{\N}})_1 = \Z$ que le morphisme $f_1$ est d\'efini par $f_1(1)$ et nous trouvons alors que l'ensemble $Hom_{Ab(\mathbf{Com} _{/\N})}(\Omega _{\N}, A)$ est aussi isomorphe au groupe ab\'elien $A_1$. 

\hfill \tf  
\end{preuve}

\vskip .5cm

\subsection{Cat\'egorie des anneaux commutatifs}

Soient $k$ un anneau commutatif et $\mathbf{Alg} _k$ la cat\'egorie des $k$-alg\`ebres commutatives, nous rappelons les d\'efinitions et les r\'esulats suivants. 

Soient $A$ une $k$-alg\`ebre et $M$ un $A$-module, alors une \emph{$k$-d\'erivation de $A$ \`a valeurs dans $M$} est un morphisme $\xymatrix @C=8mm{d : A \ar[r] & M}$ de $k$-modules v\'erifiant la r\`egle de Leibniz $d(aa')=ad(a')+a'd(a)$. 
Nous notons $Der_k(A,M)$ l'ensemble des $k$-d\'erivations, cet ensemble est muni d'une structure naturelle de $A$-module, et le foncteur 
$$\begin{array}{cccc}
Der _k (A,-) : & \mathbf{Mod}_A &  \longrightarrow & \mathbf{Ens} \\
 & M &\mapsto & Der_k(A,M),
\end{array}$$
est repr\'esent\'e par un $A$-module $\Omega _{A/k}$, appel\'e le \emph{module des diff\'erentielles de K\"ahler}, c'est-\`a-dire que nous avons un isomorphisme naturel en $M$ entre les ensembles $Hom _A(\Omega _{A/k},M)$ et $Der_k(A,M)$.  
De plus le module des diff\'erentielles $\Omega _{A/k}$ est d\'efini comme le $A$-module $I/I^2$, o\`u nous notons $I$ le noyau de la \emph{multiplication} $\xymatrix @C=8mm{m : A \otimes _k A \ar[r] & A}$. 

Pour tout $A$-module $M$ nous munissons le $A$-module $B=A \oplus M$ d'une structure de $k$-alg\`ebre appel\'ee \emph{extension de carr\'e nul}, avec un morphisme naturel $\xymatrix @C=8mm{u_B : B \ar[r] & A}$, et il existe une bijection naturelle entre l'ensemble $Der_k(A,M)$ des $k$-d\'erivations de $A$ \`a valeurs dans $M$ et l'ensemble $Hom _{{\mathbf{Alg}_k}_{/A}}(A,B)$ des morphismes de $k$-alg\`ebres de $A$ dans $B$ au-dessus de $A$. 

Nous d\'eduisons de ce qui pr\'ec\`ede le r\'esultat suivant. 

\begin{proposition} 
Pour toute $k$-alg\`ebre $A$ la cat\'egorie $Ab({\mathbf{Alg} _k}_{/A})$ des objets en groupe ab\'elien de ${\mathbf{Alg} _k}_{/A}$ est \'equivalente \`a la cat\'egorie $\mathbf{Mod}_A$ des $A$-modules. 
\end{proposition} 

\begin{preuve} 
Gr\^ace \`a la proposition \ref{prop:loi-de-groupe} appliqu\'ee \`a la loi de groupe commutatif sur un anneau nous trouvons que pour tout $B$ dans $Ab({\mathbf{Alg} _k}_{/A})$ la loi de groupe ab\'elien $m_B$ est d\'efinie par $m_B(b_1,b_2) = b_1 + b_2 - e_B(a)$ avec $a = u_B(b_1) = u_B(b_2)$. 

En appliquant l'\'egalit\'e ``$\mu _B \circ c _B= \mu _B  \circ ( \mu _B \times \mu _B) \circ \psi$" de la d\'emonstration de la proposition  \ref{prop:loi-de-groupe}  pour $\mu _B$ la multiplication dans $B$ nous en d\'eduisons que $B$ est isomorphe \`a la $k$-alg\` ebre $A \oplus M$ avec la structure de $k$-alg\`ebre \emph{extension de carr\'e nul} et $M = Ker(\xymatrix @C=8mm{u_B : B \ar[r] & A})$. 

L'\'equivalence de cat\'egories est alors une cons\'equence du corollaire \ref{cor:sous-categorie-pleine}.
 
\hfill \tf 
\end{preuve} 

\vskip .2cm 

Nous retrouvons la d\'efinition classique du module cotangent ou module des diff\'erentielles de K\" ahler $\Omega _A$ au dessus d'un anneau $A$, et les th\'eor\`emes \ref{th:premiere-suite-exacte}, \ref{th:deuxieme-suite-exacte} et \ref{th:troisieme-suite-exacte} sont des g\'en\'eralisations des r\'esultats classiques de le th\'eorie des modules des diff\'erentielles de K\"ahler. 

Le th\'eor\`eme \ref{th:premiere-suite-exacte} est une g\'en\'eralisation de la premi\`ere suite exacte fondamentale, r\'esultat qui \'enonce que pour tout morphisme $\xymatrix @C=8mm{f : A \ar[r] & B}$ dans $\mathbf{Alg} _k$ nous avons une suite exacte de $B$-modules 
$$\xymatrix @C=4mm{\Omega _{A/k} \otimes _A B  \ar@<-1pt>[rrr]^-{\delta _f} &&& \Omega _{B/k} \ar@<-1pt>[rrr] &&& \Omega _{B/A} \ar[rr] && 0 \ . }$$ 

Le th\'eor\`eme \ref{th:deuxieme-suite-exacte} est une g\'en\'eralisation du r\'esultat suivant concernant les \'epimorphismes $\xymatrix @C=8mm{f : A \ar[r] & B}$ dans $\mathbf{Alg}_k$.
La classe des \'epimorphismes est engendr\'ee par la classe des surjections $\xymatrix @C=8mm{f : A \ar[r] & B}$ avec $B=A/I$ o\`u l'id\'eal $I$ de $A$ est le noyau de $f$, et la classe des localisations $\xymatrix @C=8mm{f : A \ar[r] & B}$ avec $B=S^{-1}A$ o\`u $S$ est une partie multiplicative de $A$. 

Alors la \emph{deuxi\`eme suite exacte fondamentale} \'enonce que pour toute surjection $\xymatrix @C=8mm{f : A \ar[r] & B=A/I}$ nous avons une suite exacte de $B$-modules 
$$\xymatrix @C=4mm{I / I^2  \ar@<-1pt>[rrr] &&& \Omega _{A/k} \otimes _A B \ar@<-1pt>[rrr] ^-{\delta _f}  &&& \Omega _{B/k} \ar[rr] && 0 \ , }$$  
et la compatibilit\'e des modules des diff\'erentielles avec la localisation \'enonce que pour toute localisation $\xymatrix @C=8mm{f : A \ar[r] & B=S^{-1}A}$ nous avons un isomorphisme de $B$-modules 
$$\xymatrix @C=4mm{\delta _f : \Omega _{A/k} \otimes _A B \ar@<-1pt>[rrr] ^-{\simeq}  &&& \Omega _{B/k} \ . }$$  
Nous en d\'eduisons que pour tout \'epimorphisme $\xymatrix @C=8mm{f : A \ar[r] & B}$ l'application $\xymatrix @C=8mm{\delta _f : \Omega _{A/k} \otimes _A B \ar[r] & \Omega _{B/k}}$ est un \'epimorphisme dans $\mathbf{Mod}_B$.

Enfin le th\'eor\`eme \ref{th:troisieme-suite-exacte} est la g\'en\'eralisation du fait que les modules de diff\'erentielles commutent avec \emph{l'extension des scalaires}, 
soient $A$ une $k$-alg\`ebre, $\xymatrix @C=8mm{f : k \ar[r] & k'}$ un morphisme d'anneaux et $A' = A \otimes _k k'$, alors nous avons l'isomorphisme de $A'$-modules 
$$\xymatrix @C=4mm{\Omega _{A/k} \otimes _A A' \ar@<-1pt>[rrr] ^-{\simeq}  &&& \Omega _{A'/k'} \ . }$$  

\vskip .2cm 

Pour finir nous pouvons faire la remarque suivante. 

\begin{remark} 
Comme la cat\'egorie $Ab({\mathbf{Alg} _k}_{/A})$ est \'equivalente \`a la cat\'egorie $\mathbf{Mod}_A$, elle est ind\'ependante de l'anneau de base $k$. C'est une traduction du $(2)$ de la proposition \ref{prop:fidelite-et-equivalence}. 
\end{remark} 

\vskip 2cm

		  

\end{document}